\documentclass[11pt,reqno]{amsart}

\usepackage{etoolbox}

\newbool{bForSubmission}
\booltrue{bForSubmission}

\newbool{bIncludeStatusUpdates}
\newbool{bExperimental}
\newbool{bForUs}
\usepackage{amscd,amsmath}
\usepackage{mathrsfs}
\usepackage{amsfonts}

\usepackage{amssymb, bm, xspace}
\usepackage{enumerate}
\usepackage{setspace}

\usepackage[OT2,OT1]{fontenc}

\usepackage{datetime}

\usepackage{enumitem}

\usepackage[dvipsnames]{xcolor}

\usepackage[pagebackref=true, colorlinks=true, citecolor=blue]{hyperref}

\usepackage[margin=1.2in]{geometry}

\usepackage{mathtools}

\mathtoolsset{centercolon}

\usepackage{cleveref}

\usepackage{stackengine}
\stackMath

\usepackage{tcolorbox}
\usepackage{tikz}
\tcbuselibrary{skins}
\usetikzlibrary{shadings}
\tcbset{
    myimage/.style={
        enhanced,
        overlay={
            \begin{scope}[shift={([xshift=1mm, yshift=7mm]frame.north west)}]
            \end{scope}}}}

\tcbset{
    skin=enhanced,
    fonttitle=\bfseries,
    interior style={white},
    segmentation style={black,solid,opacity=0.2,line width=1pt}}

\newtcolorbox{TitledBox}[2][]{
    myimage,              % image 
    coltitle=black,       % title box text color
    colbacktitle=white,   % title box background color
    title=My title,
    attach boxed title to top center={
        yshift=-3mm,
        yshifttext=-1mm},
    attach boxed title to top left={
        xshift=1cm,
        yshift=-2mm},
    boxed title style={
        size=small},
    title={#2},#1}

\usepackage{pst-grad} % For gradients
\usepackage{pst-plot} % For axes
\newcommand{\TwoColumn}[6]%
{%
\begin{minipage}[#3]{#1\textwidth}%
#5%
\end{minipage}%
\begin{minipage}[#4]{#2\textwidth}%
#6%
\end{minipage}%
}

\newcommand{\TwoColumnTop}[4]%
{%
\begin{minipage}[t]{#1\textwidth}
#3
\end{minipage}%
\begin{minipage}[t]{#2\textwidth}
#4
\end{minipage}%
}

\newcommand{\PullMarginsIn}[1]%
{%
\begin{minipage}[c]{0.27\textwidth}%
\phantom{x}%
\end{minipage}%
\begin{minipage}[c]{0.46\textwidth}%
\begin{center}
#1%
\end{center}
\end{minipage}%
\begin{minipage}[c]{0.27\textwidth}%
\phantom{x}%
\end{minipage}%
}

\newbool{HaveBBM}
\booltrue{HaveBBM}

\ifbool{HaveBBM}{
	\usepackage{bbm}
    }
    {
    }

\newbool{arXivFormat}
\boolfalse{arXivFormat}
    
\newcommand{\IfarXivElse}[2]{
    \ifbool{arXivFormat}
        {#1}{#2}
    }

\renewcommand{\mathbf}[1]{\bm{#1} \textbf{ *** Use bm instead of mathbf ***}}

\newcommand{\eqn}{\begin{eqnarray}}
\newcommand{\een}{\end{eqnarray}}

\newtheorem{theorem}{Theorem}[section]
\newtheorem*{theorem*}{Theorem}				% unnumbered theorem
\newtheorem{prop}[theorem]{Proposition}
\newtheorem{lemma}[theorem]{Lemma}
\newtheorem*{lemma*}{Lemma}
\newtheorem{cor}[theorem]{Corollary}

\newtheorem{definition}[theorem]{Definition}
\newtheorem{remark}[theorem]{Remark}
\newtheorem*{remark*}{Remark}

\newtheorem{assumption}[theorem]{Assumption}
\numberwithin{equation}{section}

\newcommand{\BoldTau}
    {{\mbox{\boldmath $\tau$}}}
    
\newcommand{\bomega}
    {{\mbox{\boldmath $\omega$}}}
\newcommand{\bOmega}
    {{\mbox{\boldmath $\Omega$}}}
\newcommand{\bmu}
    {{\mbox{\boldmath $\mu$}}}

\newcommand{\BB}[1]{\ensuremath{\mathbb{#1}}}
\newcommand{\R}{{\ensuremath{\BB{R}}}}

\newcommand{\iny}{\ensuremath{\infty}}
\newcommand{\grad}{\ensuremath{\nabla}}

\newcommand{\CharFunc}{
    \ifbool{HaveBBM}{
        \ensuremath{\mathbbm{1}}
        }
        {
        \ensuremath{\bm{1}}
        }
    }

\DeclareMathOperator{\dv}{div} %
\DeclareMathOperator{\curl}{curl} %

\DeclareMathOperator{\supp}{supp} %
\DeclareMathOperator{\cond}{cond} %
\newcommand{\prt}{\ensuremath{\partial}}
\newcommand{\brac}[1]{\ensuremath{\left[ #1 \right]}}

\newcommand{\pr}[1]{\ensuremath{\left( #1 \right)}}

\DeclarePairedDelimiter{\set}{\{}{\}}

\DeclarePairedDelimiterX{\norm}[1]{\lVert}{\rVert}{#1}

\DeclarePairedDelimiterX{\abs}[1]{\lvert}{\rvert}{#1}

\newcommand\tenq[2][1]{%
	\def\useanchorwidth{T}%
	\ifnum#1>1%
		\stackunder[0pt]{\tenq[\numexpr#1-1\relax]{#2}}{\scriptscriptstyle\sim}%
	\else%
		\stackunder[1pt]{#2}{\scriptscriptstyle\sim}%
	\fi%
	}

\newcommand{\n}{{\bm{n}}}

\newcommand{\olUPlus}{\ol{U}_+ \setminus \set{0} \times \Gamma_+}
\newcommand{\onS}{\ensuremath{\text{ on } S}\xspace}
\newcommand{\onUPlus}{\ensuremath{\text{ on } \olUPlus}\xspace}

\renewcommand{\epsilon}{\varepsilon}

\newcommand{\Cal}[1]{\ensuremath{\mathcal{#1}}}
\newcommand{\al}{\ensuremath{\alpha}}
\newcommand{\la}{\ensuremath{\lambda}}

\newcommand{\diff}[2]{\frac{ d#1}{d#2}}

\newcommand{\ol}{\overline}

\renewcommand{\le}{\leqslant}
\renewcommand{\ge}{\geqslant}

\newcommand{\Holder}
    {H\"{o}lder\xspace}

\newcommand{\Gronwalls}
    {Gr\"{o}nwall's\xspace}
    
\newcommand{\NonlinearPaper}[1]{#1 of \cite{Paper2}\xspace}

\newcommand{\NoteNoTitle}[1]{
\begin{TitledBox}{}
	#1
\end{TitledBox}
}

\newcommand{\Experimental}[1]%
{%
\ifbool{bExperimental}% 
	{\bigskip
	\noindent
	\textbf{\color{Brown}*** Start experimental}\\
	#1
}
{
}
}

\definecolor{Correction}{named}{red}

\crefname{cor}{Corollary}{Corollaries} % change the second arg to \w (but
\crefname{lemma}{Lemma}{Lemmas}	       % same issue with \w, etc.

\crefname{section}{Section}{Sections}
\Crefname{section}{Section}{Sections}

\crefname{appendix}{Appendix}{Appendices}
\Crefname{appendix}{Appendix}{Appendices}

\crefname{theorem}{Theorem}{Theorems}
\Crefname{theorem}{Theorem}{Theorems}

\crefname{prop}{Proposition}{Propositions}
\Crefname{prop}{Proposition}{Propositions}

\crefname{conj}{Conjecture}{Conjectures}
\Crefname{conj}{Conjecture}{Conjectures}

\crefname{definition}{Definition}{Definitions}
\Crefname{definition}{Definition}{Definitions}

\crefname{remark}{Remark}{Remarks}
\Crefname{remark}{Remark}{Remarks}

\crefname{assumption}{Assumption}{Assumptions}
\Crefname{assumption}{Assumption}{Assumptions}

\crefformat{equation}{(#2#1#3)}
\crefrangeformat{equation}{(#3#1#4) through (#5#2#6)}
\crefmultiformat{equation}
    {(#2#1#3)}%
    { and~(#2#1#3)}
    {, (#2#1#3)}
    { and~(#2#1#3)}

\newcommand{\f}{\bm{\mathrm{f}}}

\newcommand{\g}{{\bm{\mathrm{g}}}}
\newcommand{\G}{\bm{\mathrm{G}}}
\newcommand{\h}{\bm{\mathrm{h}}}
\renewcommand{\H}{{\bm{\mathrm{H}}}}

\newcommand{\jj}{\bm{\mathrm{j}}}
\newcommand{\J}{\bm{\mathrm{J}}}

\newcommand{\uu}{{\bm{\mathrm{u}}}}   % \u is defined in some package as a check-inducing accent

\newcommand{\vv}{{\bm{\mathrm{v}}}}

\newcommand{\ww}{{\bm{\mathrm{w}}}}
\newcommand{\x}{{\bm{\mathrm{x}}}}
\newcommand{\y}{{\bm{\mathrm{y}}}}

\newcommand{\W}{{\bm{\mathrm{W}}}}
\newcommand{\X}{{\bm{\mathrm{X}}}}
\newcommand{\Y}{{\bm{\mathrm{Y}}}}
\newcommand{\z}{{\bm{\mathrm{z}}}}
\newcommand{\zz}{{\bm{\mathrm{z}}}}
\newcommand{\ZZ}{{\bm{\mathrm{Z}}}}

\newcommand{\VV}{{\bm{\Cal{V}}}}

\newcommand{\Our}{\cref{e:OurLinear}\xspace}
\newcommand{\OurOne}{\cref{e:OurLinear}$_1$\xspace}

\newcommand{\uSolSpaceN}{{C_\sigma^{N + 1, \al}(Q)}}

\newcommand{\uBoundarySpaceN}[1]{{C_{\sigma, #1}^{N + 1, \al}(Q)}}
\newcommand{\uBoundarySpaceOne}[1]{{C_{\sigma, #1}^{2, \al}(Q)}}
\newcommand{\uInputBoundarySpace}[1]{{\mathring{C}_{\sigma, #1}^{1, \al}(Q)}}
\newcommand{\uInputBoundarySpaceN}[1]{{\mathring{C}_{\sigma, #1}^{N + 1, \al}(Q)}}

\newcommand{\uInitSpaceN}{{C_\sigma^{N + 1, \al}(\Omega)}}

\newcommand{\uInputSpace}{{\mathring{C}_\sigma^{1, \al}(Q)}}
\newcommand{\uInputSpaceN}{{\mathring{C}_\sigma^{N + 1, \al}(Q)}}

\newcommand{\uInputSpaceHomoN}{\uInputBoundarySpaceN{0}}

\newcommand{\uInputSpaceArg}[1]{{\mathring{C}_\sigma^{#1 + 1, \al}(Q)}}
\newcommand{\uInputSpaceOne}{{\mathring{C}_\sigma^{2, \al}(Q)}}

\newcommand{\UnspaceN}{{C^{N + 1, \al}([0, T] \times \Gamma)}}

\newcommand{\vortSpace}{{C^\al(Q)}}
\newcommand{\vortSpaceN}{{C^{N, \al}(Q)}}
\newcommand{\vortVSpaceN}{{V^{N, \al}_\sigma(Q)}}

\newcommand{\vortInitSpace}{{C^\al(\Omega)}}
\newcommand{\vortInitSpaceN}{{C^{N, \al}(\Omega)}}
\newcommand{\vortInitSpaceOne}{{C^{1, \al}(\Omega)}}

\newcommand{\fSpaceN}{{\mathring{C}^{N + 1, \al}(Q) \cap C([0, T]; H_0)}}

\renewcommand{\time}{{\tau}}
\newcommand{\pos}{{\mbox{\boldmath $\gamma$}}}

\allowdisplaybreaks

\begin{document}
\newdateformat{mydate}{\THEDAY~\monthname~\THEYEAR}

\title
	[Linearized {E}uler with inflow, outflow]
	{The linearized 3{D} {E}uler equations with inflow, outflow}

\author[G.-M. Gie, J. Kelliher, and A. Mazzucato]
{Gung-Min Gie$^{1}$, James P. Kelliher$^{2}$, and Anna L. Mazzucato$^3$}
\address{$^1$ Department of Mathematics, University of Louisville, Louisville, KY 40292}
\address{$^2$ Department of Mathematics, University of California, Riverside, 900 University Ave., Riverside, CA 92521}
\address{$^3$ Department of Mathematics, Penn State University, University Park, PA 16802}
\email{gungmin.gie@louisville.edu}
\email{kelliher@math.ucr.edu}
\email{alm24@psu.edu}

\begin{abstract}%{{{
%	The 3D incompressible Euler equations in a bounded domain are most often
%	supplemented with impermeable boundary conditions,
%	which constrain the fluid to neither enter nor leave the domain.
	In 1983, Antontsev, Kazhikhov, and Monakhov
	published a proof of the existence and uniqueness of solutions
	to the 3D Euler equations in which on certain inflow boundary components fluid is
	forced into the domain while on other outflow components fluid is drawn out of
	the domain. A key tool they used was the linearized Euler equations
	in vorticity form.
	We extend their result on the linearized problem to multiply connected domains
	and establish compatibility conditions on the initial data that allow higher
	regularity solutions.
	% addressing an open issue in the literature.
\end{abstract}%}}}

\maketitle

%--- To add the serial comma to the list of authors in the running head;
%--- otherwise, this is exactly as generated by the amsart class
\markleft{G.-M. GIE, J. KELLIHER, AND A. MAZZUCATO}

\vspace{-2.5em}

\begin{center}
% \currfilename
Compiled on {\dayofweekname{\day}{\month}{\year} \mydate\today} at \currenttime
		
\end{center}

\bigskip

\vspace{-0.8em}

% \fbox{
\centerline{
\begin{minipage}{0.8\textwidth}
\begin{center}
\footnotesize
\renewcommand\contentsname{}
\setcounter{tocdepth}{1}		% include sections, but not subsections or below
\tableofcontents
\end{center}
% \centering
\end{minipage}
}

\section{Introduction}\label{S:Introduction}

\noindent Motivating this work and that of \cite{Paper2} is the goal of obtaining higher regularity solutions to the 3D Euler equations when fluid enters the domain through certain boundary components and exits through others---so-called inflow, outflow or injection, suction boundary conditions. Letting $\Omega$ be a bounded domain in $\R^3$ and defining the time-space domain,
\begin{align*}
	Q := (0, T) \times \Omega \text{ for a fixed but arbitrary } T > 0,
\end{align*}
we can write these equations in the form,
\begin{align}\label{e:EulerInflowOutflow}
	\begin{cases}
		\prt_t \uu + \uu \cdot \grad \uu + \grad p = \f
			& \text{in } Q, \\
		\dv \uu = 0
			& \text{in } Q, \\
		\uu(0) = \uu_0
			& \text{in } \Omega, \\		
		\uu \cdot \n = U^\n
			& \text{on } [0, T] \times \Gamma, \\
		\uu^\BoldTau  = \h
			& \text{on } [0, T] \times \Gamma_+.
	\end{cases}
\end{align}
Here $\Gamma$ is the boundary of $\Omega$; $\Gamma_+$ is the portion of the boundary on which inflow occurs; $\n$ is the outward unit normal vector; $U^\n < 0$ and $\h$ are prescribed boundary values; $\uu^\BoldTau$ is the tangential component of $\uu$; $\uu_0$ is the initial velocity; $\f$ is the external forcing. 

All proofs of existence of solutions to the Euler equations use some kind of approximation, encoded  as a sequence or as the fixed point of an operator. As the basis for one such approximation, we study the linear problem,
\begin{align}\label{e:OurLinear}
	\begin{cases}
		\prt_t \Y + \uu \cdot \grad \Y - \Y \cdot \grad \uu = \g
			&\text{in } Q, \\
		\Y = \H
			&\text{on } [0, T] \times \Gamma_+, \\
		\Y(0) = \Y_0
			&\text{on } \Omega.
	\end{cases}
\end{align}
In \Our, $\uu$ is given on $Q$, as are the initial value $\Y_0$ on $\Omega$ and the value $\H$ of $\Y$ on the inflow boundary. Should it happen that $\Y = \curl \uu$ and $\H$ is the value of the vorticity generated by the Euler equations on the inflow boundary, then $\bomega := \Y$ would be the vorticity for a solution to the Euler equations, and \OurOne would become the vorticity equation,
\begin{align}\label{e:Linear}
	\prt_t \bomega + \uu \cdot \grad \bomega - \bomega \cdot \grad \uu = \g := \curl \f.
\end{align}
We see, then, that \Our is a linearization of the vorticity formulation for \cref{e:EulerInflowOutflow}.

Employing $C^\al(Q)$ solutions to \Our, well-posedeness of \cref{e:EulerInflowOutflow} on simply connected domains for $\bomega \in C^\al(Q)$, $\al \in (0, 1)$, was obtained in Chapter 4 of \cite{AKM}. 
Here, we obtain, for any $N \ge 0$, a $C^{N, \al}(Q)^3$ solution to \Our on a multiply connected domain, which we use in \cite{Paper2} to obtain a solution to \cref{e:EulerInflowOutflow} for vorticity in $C^{N, \al}(Q)^3$ for any $N \ge 0$ on a multiply connected domain. Both the linearized \Our and the nonlinear \cref{e:EulerInflowOutflow} require
suitable compatibility conditions on the initial data to allow regularity of the solution.

Therefore, although \cref{e:EulerInflowOutflow} motivates our work, we restrict our attention to \Our.

\subsection*{The key difficulty} If $\g \equiv 0$, \OurOne shows that $\Y$ is the pushforward (transport with stretching) of $\Y_0$ (we explain this in detail in \cref{S:Pushforward}). The trajectories of the flow map for $\uu$ play a central role. If we assume that $\uu \cdot \n = 0$ on the boundary then this is nearly the complete story other than accounting for forcing, and \Our is solved in an entirely classical manner.

With inflow of fluid from the boundary, however, we must insure that the values of $\Y = \H$ coming from the inflow meet seamlessly enough with those coming from the initial data $\Y_0$ so that the desired regularity of the solution is obtained. This is the primary complication we face in solving \Our.

\subsection*{Inflow, outflow}
We assume that $\Gamma := \prt \Omega$ has at least $C^2$ regularity and has a finite number of components, with $\Gamma_+$, $\Gamma_-$, $\Gamma_0$ a partition of the boundary components into those on which inflow, outflow, no-penetration boundary conditions hold, respectively. That is, defining
\begin{align*}
	U^\n
		:= \uu \cdot \n,
\end{align*}
we require that
\begin{align}\label{e:UnThree}
	U^\n < 0 \text{ on } [0, T] \times \Gamma_+, \qquad
	U^\n > 0 \text{ on } [0, T] \times \Gamma_-, \qquad
	U^\n = 0 \text{ on } [0, T] \times \Gamma_0.
\end{align}
(We allow $\Gamma_0 = \emptyset$ or $\Gamma_0 = \Gamma$---see \cref{R:NoPenetrationBCs}.)
Moreover, $\dv \uu = 0$ imposes the constraint,
\begin{align}\label{e:UnConstraint}
	\int_{\Gamma_+} U^\n = - \int_{\Gamma_-} U^\n.
\end{align}

\PullMarginsIn{%
\NoteNoTitle
{
	Throughout, we fix $\al \in (0, 1)$.
}%
}

\subsection*{Linear problem as a tool}
With sufficient regularity we interpret \Our classically, but for our lowest regularity solutions \OurOne must be treated weakly, as equality in the sense of distributions on $Q$. In all cases, we will construct and treat $\Y$ as a Lagrangian solution to \Our, though to even define what we mean by such solutions will require the development of some technology because of the inflow of $\H$ from $\Gamma_+$ (this leads ultimately to \cref{D:LagrangianSolution}).

We can view \cref{e:Linear} as a special case of \cref{e:OurLinear}, in which $\Y = \bomega = \curl \uu$ and $\H$ is derived from the pressure, as done in \cite{AKM, Paper2}. In \cref{e:Linear}, $\bomega$ is a curl and so, in particular, is divergence-free, whereas this is not assumed for $\Y$ in \Our. We will show, however, that if $\dv \Y_0 = 0$ and $\H$ satisfies the condition in \cref{e:RangeOfCurlCond}, obtained formally be restricting \OurOne to $\Gamma_+$, then $\dv \Y(t)$ will be zero for $t \in [0, T]$. Moreover, we show that if $\Y_0$ is in the range of the curl then $\Y(t)$ remains in the range of the curl for $t \in [0, T]$, which requires additional work only because $\Gamma$ has multiple components (unless, perhaps, $\Gamma = \Gamma_0$).

The analysis of \cref{e:OurLinear} in \cite{AKM} focused on $C^\al(Q)$ regularity. We are concerned with obtaining $C^{N, \al}(Q)$ regularity of solutions to both the linear (in this paper) and nonlinear (in \cite{Paper2}) problems for any integer $N \ge 0$. To accomplish this, we must discover the right compatibility conditions on the initial data and on $\H$. For $N = 0$ the conditions, obtained in \cite{AKM}, are simply that $\H(0) = \Y_0$ on $\Gamma_+$. We will show that these conditions have a natural generalization to all $N \ge 0$, most cleanly stated below in the necessary and sufficient conditions in \cref{e:CondNAlt}.

\subsection*{Some function spaces}
Let $V$ be an open subset of $\R^d$, $d \ge 1$. We define the classical \Holder space $C^\al(V)$ to be all measurable real-valued functions on $V$ for which
\begin{align*}
	\norm{f}_{C^\al(\Omega)}
		:= \norm{f}_{L^\iny(\Omega)} + \norm{f}_{\dot{C}^\al(\Omega)} < \iny, \quad
	\norm{f}_{\dot{C}^\al(V)}
		:= \sup_{x \ne y \in V}
			\frac{\abs{f(x) - f(y)}}{\abs{x - y}^\al}.
\end{align*}
For any integer $N \ge 0$ we define the Banach space $C^{N, \al}(V)$ with the norm
\begin{align*}
	\norm{f}_{C^{N, \al}(V)}
		&:= \sum_{\abs{\gamma} \le N} \norm{D^\gamma f}_{L^\iny}
			+ \sum_{\abs{\gamma} = N} \norm{D^\gamma f}_{\dot{C}^\al(V)}.
\end{align*}
We also allow $f$ to be vector- or matrix-valued, but will not make a notational distinction.

For the time-space domain $Q$, we define
\begin{align*}
	\mathring{C}^{N + 1, \al}(Q)
		&:= \set{\vv \colon Q \to \R \colon
				\prt_t^j D_\x^\gamma \vv \in C^\al(Q)^3,
				j + \abs{\gamma} \le N + 1, j \le N},
\end{align*}
endowed with the natural norm based upon its regularity. That is, $\mathring{C}^{N + 1, \al}(Q)$ is the same as $C^{N + 1, \al}(Q)$, but with one less time than spatial derivative of regularity.

We call $\beta \in C^{N + 1, \al}([0, T] \times \Gamma)$ a \textit{proper inflow, outflow boundary value} if it satisfies the same conditions as $U^\n$ does in \cref{e:UnThree,e:UnConstraint}. We then define the affine spaces
\begin{align}\label{e:BoundarySpace}
	\begin{split}
	\uBoundarySpaceN{\beta}
		&:= \set{\vv \in C^{N + 1, \al}(Q) \colon \dv \vv = 0, \vv \cdot \n = \beta
				\text{ on } [0, T] \times \Gamma}, \\
	\uInputBoundarySpaceN{\beta}
		&:= \set{\vv \in \mathring{C}^{N + 1, \al}(Q) \colon \dv \vv = 0, \vv \cdot \n = \beta
				\text{ on } [0, T] \times \Gamma}.
	\end{split}
\end{align}

Primarily, we will utilize the following:
\begin{align}\label{e:CsigmaSpaces}
	\begin{split}
	\uInitSpaceN
		&:= \set{\vv \in C^{N + 1, \al}(\Omega) \colon \dv \vv = 0, \vv \cdot \n = U^\n
					\text{ on } \Gamma}, \\
	\uSolSpaceN
		&:= \uBoundarySpaceN{U^\n},
		\qquad
	\uInputSpaceN
		:= \uInputBoundarySpaceN{U^\n},
	\end{split}
\end{align}
where we suppose that $U^\n$ is at least as regular as $C^{N + 1, \al}([0, T] \times \Gamma)$.
Observe that only the normal component of $\vv$ is specified on the boundary.

We will also use the classical space,
\begin{align}\label{e:HSpace}
	H
		&:= \set{\vv \in L^2(\Omega)^3 \colon \dv \vv = 0, \, \vv \cdot \n = 0
				\text{ on } \Gamma}
		= H_0 \oplus H_c,
\end{align}
where the $L^2$-orthogonal subspaces $H_c$, $H_0$ of $H$ are defined by
\begin{align}\label{e:HcH0}
	H_c &:= \set{\vv \in H: \curl \vv = 0}, \quad
	H_0 := H_c^\perp.
\end{align}

\subsection*{Regularity assumptions on the data} We specify the regularity of the initial data, boundary value $U^\n$, $\H$ on inflow, forcing $\g$, and velocity field $\uu$, as follows:

\begin{definition}\label{D:NReg}
	We say the data has regularity $N$ for integer $N \ge 0$ if the following hold:
	\begin{itemize}
		\item
			$\Gamma$ is $C^{\max\set{N + 2, 3}, \al}$,
			$U^\n \in \UnspaceN$;
			
		\item
			$\g \in C^\al(Q)$ if $N = 0$,
			$\g \in \mathring{C}^{N, \al}(Q)$ if $N \ge 1$;
						
		\item
			$\Y_0 \in \vortInitSpaceN$,
			$\H \in C^{N, \al}([0, T] \times \Gamma_+)$;
			
		\item
			$\uu \in \uInputSpaceN$.		
	\end{itemize}
\end{definition}

In \cite{AKM}, $\H$ and $U^\n$ are assumed to have one more derivative of regularity than we have assumed here for data regularity $N = 0$. Higher regularity of $\H$ is required, as we will see in \cref{T:LinearExistence}, to insure that $\Y(t)$ remains in the range of the curl if $\Y_0$ is in the range of the curl. Since \cite{AKM} analyzes solutions to the Euler equations, such higher regularity is needed. Moreover, the need to properly control the pressure for the nonlinear problem, which is directly related to the production of vorticity on the boundary, requires higher regularity of both $\H$ and $U^\n$ in \cite{AKM} (and in \cite{Paper2}). But that is not an issue for the linear problem we treat here.

For most of the analysis we make, $\Gamma$ being $C^{N + 2, \al}$ is sufficient. In the proof of \cref{L:VelDense01}, however, which we apply only for data regularity $N = 0$, we need one more derivative of regularity. See \cref{R:GammaHigherRegularity}.

\subsection*{Compatibility conditions}
To obtain the regularity of solutions to \Our, we need to impose compatibility conditions. We define the conditions for $N = 0$ and $N = 1$ as
\begin{align}\label{e:Cond01}
	\begin{array}{ll}
		\cond_0: &\H(0) = \Y_0 \text{ on } \Gamma_+, \\
		\cond_1: &\cond_0 \text{ and }
			\prt_t \H|_{t = 0} + \uu_0 \cdot \grad \Y_0 - \Y_0 \cdot \grad \uu_0 - \g(0) = 0
				\text{ on } \Gamma_+,
	\end{array}
\end{align}
where $\uu_0 := \uu(0)$. We can view $\cond_1$ formally as saying that \OurOne holds at time zero on $\Gamma_+$, where $\Y_0 = \H(0)$ by $\cond_0$. Indeed, we could write $\cond_N$ for all $N \ge 1$ suggestively as
\begin{align}\label{e:CondNAlt}
	\begin{array}{ll}
		\cond_N: &\prt_t^j \H|_{t = 0} = \prt_t^j \Y|_{t = 0}
				\text{ on } \Gamma_+
				\text{ for all } 0 \le j \le N,	
	\end{array}
\end{align}
where we replace $\prt_t \Y$ by the form it would have were it a solution to \OurOne.
Or, spelled out just a little more,
\begin{align*}
	\begin{array}{ll}
		\cond_N: &\cond_{N - 1} \text{ and }
			\prt_t^N \H|_{t = 0}
				+ \prt_t^{N - 1} [\uu \cdot \grad \Y - \Y \cdot \grad \uu - \g]_{t = 0}
				= 0
				\text{ on } \Gamma_+.
	\end{array}
\end{align*}

For $N = 2$, for instance, we would first write,
\begin{align*}
	0
		&= \prt_t^2 \H|_{t = 0}
				+ \prt_t [\uu \cdot \grad \Y - \Y \cdot \grad \uu - \g]_{t = 0} \\
		&= \prt_t^2 \H(0)
				+ \prt_t \uu(0) \cdot \grad \Y_0 + \uu_0 \cdot \grad \prt_t \Y(0)
				- \prt_t \Y(0) \cdot \grad \uu_0 - \Y_0 \cdot \grad \prt_t \uu(0)
				- \prt_t \g(0),
\end{align*}
then replace each $\prt_t \Y(0)$ with $\g(0) - \uu_0 \cdot \grad \Y_0 + \Y_0 \cdot \grad \uu_0$, the value it would have were \OurOne to hold.

\subsection*{Types of solutions}
We will be concerned with both Lagrangian solutions (\cref{D:LagrangianSolution}) and with Eulerian solutions, classical as well as weak, as in \cref{D:Weak}.
In this definition, $\Y$ has more than sufficient time and boundary regularity to avoid the need to enforce the initial and boundary conditions weakly. With sufficient regularity, $\dv(\Y \otimes \uu) =\uu \cdot \grad \Y$ (using $\dv \uu = 0$). Since we only assume $\Y \in C^\al(Q)$, 
$\dv(\Y \otimes \uu)$ is defined in the sense of distributions, given that $\Y \otimes \uu$ is an integrable function.

\begin{definition}\label{D:Weak}
	We say that $\Y \in \vortSpace$ is a weak (Eulerian) solution to \Our if
	$\Y = \H$ on $[0, T] \times \Gamma_+$, $\Y(0) = \Y_0$, and
	$\prt_t \Y + \dv(\Y \otimes \uu) - \Y \cdot \grad \uu = \g$
	in $\Cal{D'}(Q)$.
\end{definition}

In short, we will find that the Lagrangian solution is a weak Eulerian solution if $\Y$ is in the range of the curl, and weak Eulerian solutions lying in the range of the curl are unique.
We have a particular concern over the Lagrangian versus Eulerian solution because we use both formulations when we apply our results in \cite{Paper2}---specifically, when $\Y$ is the vorticity of some vector field. The vast majority of the estimates come from the Lagrangian formulation, but the velocity formulation, which is also needed, is recovered from the weak Eulerian formulation, not from the Lagrangian, so they need to be the same solution.

\subsection*{Main results}
Our main results are \cref{T:LinearExistence,T:VelocityForm}.

\begin{theorem}\label{T:LinearExistence}
	Assume that the data has regularity $N$ for some $N \ge 0$ and that
	$\cond_N$ holds.
	There exists
	a solution $\Y$ to \Our in $C^{N, \al}(Q)$ with
	\begin{align}\label{e:RegBound}
		\norm{\Y}_{C^{k, \al}(Q)}
			&\le C(T) \brac{\norm{\uu}_{\uInputSpaceArg{k}}
				+ \norm{\H}_{C^{k, \al}([0, T] \times \Gamma_+)}
				+ \norm{\g}_{C^{\max\set{k - 1, 0}, \al}(Q)} T}
	\end{align}
	for all $0 \le k \le N$, where $C(T)$ also depends upon
	$\norm{\Y_0}_{C^{k, \al}(\Omega)}$.
	The solution $\Y$ is Lagrangian and weak Eulerian for
	all $N \ge 0$. For $N \ge 1$, the solution is also the unique classical Eulerian
	solution.

	Moreover, suppose that $\Y_0$ and $\g(t)$ for all $t$ are in the
	range of the curl,
	\begin{align}\label{e:HigherH}
		\H \in C^{\max\set{N, 1}, \al}([0, T] \times \Gamma_+),
	\end{align}
	and (our notation is defined in \cref{e:VectorComponentsNotation})
	\begin{align}\label{e:RangeOfCurlCond}
	\prt_t H^\n
		+ \dv_\Gamma[H^\n \uu^\BoldTau  - U^\n \H^\BoldTau]
		- \g \cdot \n
		= 0
			\text{ on } (0, T] \times \Gamma_+.
	\end{align}
	Then $\Y(t)$ remains in the range of the curl for all $t \in [0, T]$;
	also, for $N = 0$, $\Y$ is the unique weak Eulerian solution.
\end{theorem}

We make a few comments on \cref{T:LinearExistence}:
\begin{itemize}
\item
In \cref{e:RangeOfCurlCond}, $\dv_\Gamma$ is the divergence-operator along $\Gamma_+$ (see \cref{S:Vorticity}).

\item
We give the definition of a Lagrangian solution in \cref{D:LagrangianSolution}. It is not entirely classical because values of $\H$ are brought into the domain from the inflow boundary.

\item
Because of \cref{e:HigherH}, both $\g \cdot \n$ and $\prt_t H^\n$ are in $C^\al([0, T] \times \Gamma_+)$ in \cref{e:RangeOfCurlCond}. Hence, it is implicit in \cref{e:RangeOfCurlCond} that $\dv_\Gamma[H^\n \uu^\BoldTau  - U^\n \H^\BoldTau]$ is in $C^\al([0, T] \times \Gamma_+)$.

\item
Restricted to $t = 0$, \cref{e:RangeOfCurlCond} is the normal component of $\cond_1$, as we can see by examining the proof of \cref{P:YInRangeOfCurl}.
\end{itemize}

We also have a velocity formulation of \Our, as given in \cref{T:VelocityForm}.

\begin{theorem}\label{T:VelocityForm}
	Assume that the data has regularity $N$ for some $N \ge 0$ and that
	$\cond_N$, \cref{e:HigherH}, and \cref{e:RangeOfCurlCond} hold.
	Let $\f \in \fSpaceN$ and set
	$
 		\g := \curl \f.
	$
	Let $\bomega = \Y$ and $\bomega_0 = \Y_0$.
	There exists a unique $\vv \in \uSolSpaceN$ with $\curl \vv = \bomega$
	and mean-zero pressure field $\pi$ with $\grad \pi \in C^{N, \al}(Q)$
	for which
	\begin{align}\label{e:LinearVelEq}
		\prt_t \vv + \uu \cdot \grad \vv - \uu \cdot (\grad \vv)^T + \grad \pi
			= \f.
	\end{align}
	The harmonic component of $\vv$ is given explicitly in \cref{L:Vc}.
	
	Suppose, further, that $\bomega = \curl \uu$. Then
	setting
	$p = \pi - (1/2) \abs{\uu}^2$, $(\uu, p)$ satisfies
	\cref{e:EulerInflowOutflow}$_{1-4}$ (but not \cref{e:EulerInflowOutflow}$_5$)
	with an additional, harmonic forcing term. If, further, $\vv = \uu$ then this additional
	term does not appear and the solution is unique, even for $N = 0$.
\end{theorem}

\begin{remark}
	\cref{T:VelocityForm} shows that in velocity formulation, the solution to the linearized
	Euler equations is unique, even for $N = 0$.
	The last paragraph of \cref{T:VelocityForm} considers what would happen were the linear solution
	to be a solution to the Euler equations; it does not establish the existence of such
	a solution---that is done in \cite{Paper2}. Also, because of how $\vv$ is recovered from the
	vorticity, $\vv$ and $\uu$ have matching normal
	components on the boundary, but will not, in general, have matching tangential components.
\end{remark}

\subsection*{Prior work}
The primary reference for our work is Chapter 4 of \cite{AKM}. Although we express things differently, the material in \cref{S:FlowMap,S:Pushforward,S:Solutions,S:Vorticity,S:Velocity} clearly bears the imprint of \cite{AKM}. Our motivation in this work and \cite{Paper2} is to ultimately extend the results of Chapter 4 of \cite{AKM}, which are for simply connected domains and $N = 0$ regularity, to obtain solutions to the Euler equations with $N \ge 0$ regularity with suitable compatibility conditions.

Section 1.4 of \cite{Mamontov2010} contains an extensive survey of results, both 2D and 3D, related to the problem we are studying here. Petcu \cite{Pet06} presents a version of the argument in Chapter 4 of \cite{AKM} specialized to a 3D channel with inflow and outflow constant in space and time.

Boyer and Fabrie \cite{BoyerFabrie2006,BoyerFabrie2013} (in particular, see the final two chapters of \cite{BoyerFabrie2013}) treat an analogous setup to ours for transport without stretching, studying weak solutions. They do not restrict inflow and outflow to lie on full components of the boundary, which would seem to make classical solutions impossible to obtain. In a related vein, see also the recent works \cite{BravinSueur2021,NoisetteSueur2021}. We also mention the works \cite{CST10,GHJT18,HT09} on different linear problems.

Finally, we note that the need for the higher regularity results of \cref{T:LinearExistence,T:VelocityForm} is explicitly stated in \cite{TemamWang2002,GHT2}, where such results are used to obtain high-order expansions of solutions to the Navier-Stokes equation asymptotically in terms of the viscosity.

\subsection*{Organization of this paper}
We start by developing some necessary tools. In \cref{S:BSLaw} we describe how to recover a velocity in $C_\sigma^{N + 1, \al}(\Omega)$ from a vorticity in $C^{N, \al}(\Omega)$ and introduce the concepts we need to treat multiply connected domains. In \cref{S:FlowMap}, we develop the properties of the flow map, which we will need throughout the rest of the paper.

We develop the idea of a pushforward of a velocity field with boundary conditions on inflow, a non-classical construct, in \cref{S:Pushforward}, using it in \cref{S:Solutions} to define what we mean by a Lagrangian solution. We give the core of the proof of \cref{T:LinearExistence} in \cref{P:RegNFromS,P:LinearExistenceFirstPart}, though we defer to \cref{S:N1RegCond} a key part of it, giving an equivalent form of $\cond_1$: it is of a much different flavor than the rest of the proof and would only distract from it.

The proof that if $\Y_0$ lies in the range of the curl then $\Y(t)$ remains in that range for $t > 0$ is given in \cref{S:Vorticity}. The results in \cref{S:Solutions,S:Vorticity}, specifically \cref{P:LinearExistenceFirstPart,P:YInRangeOfCurl}, together yield \cref{T:LinearExistence}. The proof of \cref{T:VelocityForm} is presented in \cref{S:Velocity}. In the appendices, we list some standard \Holder space estimates and give details on the continuity of the Biot-Savart law, which we referenced in \cref{S:BSLaw}.

\subsection*{On notation} Our notation, while fairly standard, has a few subtleties.
If $M$ is a matrix, $\cramped{M^i_k}$ refers to the entry in row $i$, column $k$ of $M$; $\cramped{v^i}$ refers to the $\cramped{i^{th}}$ entry in the vector $\vv$, which we always treat as a column vector for purposes of multiplication. If $M$ and $N$ are the same size matrices then $\cramped{M \cdot N := M^i_k N^i_k}$, where here, as always, we use implicit summation notation. If $\uu$ and $\vv$ are vectors then the matrix $\uu \otimes \vv$ has components $[\uu \otimes \vv]^i_k := u^i v^k$. We define the divergence of a  matrix row-by-row, so $\dv M$ is the column vector with components $[\dv M]^i = \prt_k M^i_k$. Hence, $[\dv [\uu \otimes \vv]]^i = \dv [\uu \otimes \vv]^i = \prt_k (u^i v^k)$, where $\prt_k$ is the derivative with respect to the $k^{th}$ spatial variable. For any vector field $\vv$ defined on $\Gamma$, we define its normal and tangential components,
\begin{align}\label{e:VectorComponentsNotation}
	v^\n := \vv \cdot \n, \quad
	\vv^\n := v^\n \n, \quad
	\vv^\BoldTau := \vv - \vv^\n.
\end{align}

The operators $\grad$, $\dv$ are the gradient, divergence with respect to the spatial variables only. For vector fields $\uu$ and $\vv$, we will interchangeably write $\uu \cdot \grad \vv$ and $\grad \vv \, \uu$, each of which is a vector whose $i^{th}$ component is $u^k \prt_k v^i$.
When applied to a function $\eta$ that includes two time variables (as our flow maps will), we write $\prt_{t_1} \eta$, $\prt_{t_2} \eta$ to mean the derivative with respect to the first, second time variable.
We will sometimes treat time-space as a four-dimensional variable, defining the operator
\begin{align}\label{e:D}
	D := (\prt_t, \grad),
\end{align}
noting that if $\vv$ is a vector field then $D \vv$ is a $3 \times 4$ matrix field.

Finally, $H^{k, p}(\Omega)$ and $H^s(\Omega)$ are the standard $L^p$- and $L^2$-based Sobolev spaces.

\section{Recovering velocity from vorticity}\label{S:BSLaw}

\noindent In this section, we develop the operator in \cref{e:KUDef} that we will use to recover an appropriate divergence-free vector field from its vorticity (curl), with its regularity given in  \cref{L:KRegular}. Other than in the proof of one technical, but important, result that we will defer to \cref{A:KTechnical}, we will need only a few facts regarding the Hodge decomposition for multiply connected domains, which we now summarize. 

Assume that $\Omega$ is connected and that $\Gamma$ is at least $C^2$-regular. Let $\Gamma_1, \ldots, \Gamma_{b + 1}$, be the $b + 1$ components of $\Gamma$ with $\Gamma_{b + 1}$ the boundary of the unbounded component of $\Omega^C$.
We define the \textit{external flux} of $\bomega$ through $\Gamma_i$ as
\begin{align}\label{e:ExternalFlux}
	\Phi^\Gamma_i(\bomega) := \int_{\Gamma_i} \bomega \cdot \n.
\end{align}

As we discuss in \cref{A:KTechnical}, the space $H_0$ consists of those elements of $H$ having vanishing \textit{internal} fluxes, a characterization first given by Helmholtz (see the historical comments in \cite{CantarellaDeTurckGluck2002}). For our purposes, we do not need an explicit characterization of these spaces, only the definitions of the spaces $H$, $H_c$, and $H_0$ in \cref{e:HSpace,e:HcH0} and the fact that $H_c$ is finite-dimensional. Employing elliptic regularity theory, \cref{L:TamingHc} easily follows:
\begin{lemma}\label{L:TamingHc}
	Assume that $\Gamma$ is $C^{n, \al}$-regular, $n \ge 2$, and
	let $X$ be any function space that contains $C^{n, \al}(\Omega)^3$.
	For any $\vv \in H$,
	\begin{align*}
		\norm{P_{H_c} \vv}_X \le C(X) \norm{\vv}_H
	\end{align*}
	and if also $\vv \in X$ then
	\begin{align*}
		\norm{\vv}_X \le \norm{P_{H_0} \vv}_X + C(X) \norm{\vv}_H, \quad
		\norm{P_{H_0} \vv}_X \le \norm{\vv}_X + C(X) \norm{\vv}_H.
	\end{align*}
\end{lemma}

For any $\vv \in H = H_0 \oplus H_c$, we call $P_{H_c} \vv$ the \textit{harmonic component or part of $\vv$}.
(Note that $\Delta \vv_c = 0$ for any $\vv_c \in H_c$, though unless $\Omega$ is simply connected, there are also $\vv \in H_0$ for which $\Delta \vv = 0$.)

The following is a classical trace theorem (see Theorem 1.2 p. 7 of \cite{T2001}):
\begin{lemma}\label{L:Trace}
	There is a continuous trace from the space of $L^2$ vector fields on $\Omega$ having divergence
	in $L^2(\Omega)$ to $H^{-\frac{1}{2}}(\Gamma)$.
\end{lemma}

We write $\curl H^1(\Omega)^3$ for the image of $H^1(\Omega)^3$ under the curl operator and say that a vector field is \textit{in the range of the curl} if it lies in $\curl H^1(\Omega)^3$. We can characterize it as in \cref{T:VectorPotential}.

\begin{theorem}\label{T:VectorPotential}
	We have,
	\begin{align*}
		\curl H^1(\Omega)^3
			= \curl (H_0 \cap H^1(\Omega)^3)
			= \set{\bomega \in L^2(\Omega) \colon \dv \bomega = 0,
				\Phi^\Gamma_i(\bomega) = 0 \text{ for all } i}.
	\end{align*}
	Moreover, there exists a continuous operator
	$K$ from $\curl H^1(\Omega)^3 \subseteq L^2(\Omega)^3$ to $H_0 \cap H^1(\Omega)^3$
	with the property that $\uu = K[\bomega]$ is the unique element of
	$H_0 \cap H^1(\Omega)^3$ for which $\curl \uu = \bomega$.		
\end{theorem}
\begin{proof}
	This follows from Theorems 3.5 and 3.12 of \cite{ABDG1998}.
	Observe that \cref{L:Trace} allows
	the external flux $\Phi^\Gamma_i(\bomega)$ to be defined.
\end{proof}

To recover a velocity field $\uu$ satisfying $\uu \cdot \n = U^\n$, we let $\VV = \grad \varphi$, where $\varphi$ is the unique mean-zero solution to
\begin{align*}
	\begin{cases}
		\Delta \varphi = 0
			&\text{in } \Omega, \\
		\grad \varphi \cdot \n = U^\n
			&\text{on } \Gamma.
	\end{cases}
\end{align*}
Observe that $\dv \VV = 0$, $\curl \VV = 0$, and $\VV \cdot \n = U^\n$ on $\Gamma$.
(Note that if $U^\n = 0$ then $\VV \equiv 0$.) Then we define
\begin{align}\label{e:KUDef}
	K_{U^\n}[\bomega]
		:= K[\bomega] + \VV.
\end{align}

Define the solution space for vorticity,
\begin{align*}
	\vortVSpaceN
		:= \set{\bomega \in C^{N, \al}(Q)^3 \colon
			\bomega(t) \in \curl H^1(\Omega)^3 \text{ for all } t \in [0, T]}.
\end{align*}
From \cref{T:VectorPotential}, $\bomega \in \vortVSpaceN$ is equivalent to $\bomega \in \vortSpaceN$ lying in the range of the curl.

\begin{lemma}\label{L:KRegular}
	Assume the data has regularity $N \ge 0$.
	Then $K_{U^\n}$ maps $\vortInitSpaceN \cap \curl H^1(\Omega)^3$ continuously
	onto $\uInputSpaceN \cap (H_0 + \VV(t))$ and maps $\vortVSpaceN$ continuously
	onto
	$
		\uInputSpaceN \cap C([0, T]; H_0 + \VV(t))
	$.
\end{lemma}
\begin{proof}
	Follows from \cref{e:KUDef}, the regularity of $\VV$, \cref{L:TamingHc}, and
	\cref{C:BSLawLikeBound}.
\end{proof}

We also have the Helmholtz decomposition, which we use to establish the continuity of the Leray projector. We give its proof, since the continuity of the decomposition in \Holder spaces, especially for $N = 0$, is not particularly accessible in the literature.
\begin{lemma}\label{L:LerayContinuity}
	Assume $\Gamma$ is $C^{N + 1, \al}$ for $N \ge 0$.
	Given any $\bomega \in C^{N, \al}(\Omega)^3$, we have
	$\bomega = \curl \vv + \grad q$ for some unique $\vv \in H_0 \cap C^{N + 1, \al}(\Omega)^3$
	and mean-zero $q \in C^{N + 1, \al}(\Omega)$, the maps from
	$\bomega$ to $\vv$ and $\bomega$ to $q$ being continuous in their respective spaces.
	Moreover, the Leray projector, $P_H$, is continuous as a map from $C^{N, \al}(\Omega)^3$ to
	$C^{N, \al}(\Omega)^3 \cap H$.
\end{lemma}
\begin{proof}
	Let $\bomega = P_H \bomega + \grad q$ be the Leray-Helmholtz decomposition of $\bomega$.
	By \cref{L:KRegular} (applied using $\VV \equiv 0$),
	$P_H \bomega = \curl \vv$ for some unique $\vv \in H_0 \cap C^{N + 1, \al}(\Omega)^3$,
	with $\bomega \mapsto \vv$ being a continuous map. Hence,
	$\bomega = \curl \vv + \grad q$, and the maps,
	$\bomega \mapsto \vv$, $\bomega \mapsto q$, and $\bomega \mapsto P_H \bomega$
	have the stated continuity.
\end{proof}

\begin{lemma}\label{L:VelDense01}
	Assume that the data has $N = 0$ regularity. 	
	There exists a sequence $(U_i^\n)$ of proper inflow, outflow boundary values
	and a sequence $(\uu_i)$ in $\uBoundarySpaceOne{U_i^\n}$
	converging to $\uu$ in $\mathring{C}^{1, \al}(Q)$.
	(The space $\uBoundarySpaceOne{U_i^\n}$ is defined in \cref{e:BoundarySpace}.)
\end{lemma}
\begin{proof}
	Let $(U_i^\n)$ be a sequence of vector fields mollified along
	$\Gamma$ so that $U_i^\n \in C^{2, \al}([0, T] \times \Gamma)$
	with $U_i^\n \to U^\n$ in $C^{1, \al}([0, T] \times \Gamma)$.
	This is possible,
	since we assumed that $\Gamma$ is $C^{2, \al}$. Then for all
	sufficiently large $i$, $U_i^\n$ will satisfy the same conditions as $U^\n$ does in
	\cref{e:UnThree} and, after possibly adjusting the value of
	$U_i^\n$ on one boundary component by adding a constant value $c_i$ to it,
	with $c_i \to 0$, each $U_i^\n$ will also satisfy \cref{e:UnConstraint}.
	Let $\VV_i = \grad \varphi_i$, where
	$\varphi_i$ solves
	\begin{align*}
		\begin{cases}
				\Delta \varphi_i = 0
				&\text{in } \Omega, \\
			\grad \varphi_i \cdot \n = U_i^\n
				&\text{on } \Gamma.
		\end{cases}
	\end{align*}
	Then $\VV_i \in C^{2, \al}(Q)$ with $\VV_i \to \VV$ in $C^{1, \al}(Q)$
	by elliptic regularity theory.
	
	We have $\ww := \uu - \VV \in \uInputBoundarySpace{0}$, which is the space
	$\uInputSpace$, but with
	$\ww \cdot \n = 0$ on $\Gamma$.
	Extend $\ww$ to
	$\mathring{C}^{1, \al}(\R \times \R^3)$ using
	an extension operator like that in Theorem 5', chapter VI of \cite{Stein1970}.
	Mollify $\ww$ in time and space and apply the Leray projector $P_H$,
	giving, via \cref{L:LerayContinuity} a sequence $\ww_i \in \uBoundarySpaceOne{0}$
	with $\ww_i \to \ww$ in $\uInputBoundarySpace{0}$.
	
	Finally, let $\uu_i = \ww_i + \VV_i$.
\end{proof}

\begin{remark}\label{R:GammaHigherRegularity}
	 In the proof of \cref{L:VelDense01}, we applied \cref{L:LerayContinuity} with
	 $N = 2$, which required that $\Gamma$ (through \cref{C:BSLawLikeBound}) be $C^{3, \al}$.
	 This is the only place
	 in this paper in which we required higher than a $C^{2, \al}$ boundary
	 for data regularity 0.
\end{remark}

\section{The flow map}\label{S:FlowMap}

\noindent In this section, we assume that for some $N \ge 0$ and fixed $T > 0$,
$
	\uu \in \uInputSpaceN
$.
We also assume that $\Gamma$ is $C^{N + 2, \al}$-regular. We will obtain estimates related to the flow map for $\uu$.

For convenience, we first extend $\uu$ to be defined on all of $\R \times \R^3$ using an extension operator like that in Theorem 5', chapter VI of \cite{Stein1970}. This extension need not be divergence-free. This extension will allow us to use classical results on flow maps without undue concern as to their domain and codomain, though it is only the value of $\uu$ on $\ol{Q}$ that ultimately concerns us.

Define $\eta \colon \R \times \R \times \R^3 \to \R^3$ to be the unique flow map for $\uu$, so that
\begin{align}\label{e:FlowId}
	\begin{split}
	&\prt_{t_2} \eta(t_1, t_2; \x) = \uu(t_2, \eta(t_1, t_2; \x)),  \quad \eta(t_1, t_1; \x) = \x; \\
	&\eta(t_1, t_2; \x)
		= \x + \int_{t_1}^{t_2} \uu(s, \eta(t_1, s; \x)) \, ds.
	\end{split}
\end{align}
That is, $\eta(t_1, t_2; \x)$ is the position that a particle starting at time $t_1$ at position $\x \in \R^3$ will be at time $t_2$ as it moves in the velocity field $\uu$. We allow $t_2$ to be greater than, equal to, or less than $t_1$, accounting for movement forward and backward in time.
(The properties of flow lines within $\Omega$, which are all we ultimately care about, do not depend upon the specific extension of $\uu$ we employ.) We also have
\begin{align}\label{e:prtt1eta}
	\begin{split}
	\prt_{t_1} \eta(t_1, t_2; \x)
		&= - \uu(t_1, \eta(t_1, t_1; \x))
			+ \int_{t_1}^{t_2} \prt_{t_1} \uu(s, \eta(t_1, s; \x)) \, ds \\
		&= - \uu(t_1, \x)
			+ \int_{t_1}^{t_2} \grad \uu(s, \eta(t_1, s; \x))
					\prt_{t_1} \eta(t_1, s; \x) \, ds.
	\end{split}
\end{align}
Moreover, by the very definition of the flow map (and its uniqueness),
\begin{align}\label{e:GroupComp}
	\eta(t_2, t_3; \eta(t_1, t_2; \x))
		= \eta(t_1, t_3; \x).
\end{align}

\cref{L:etaFlow} shows that for fixed $t_2$ and $\x$, $\eta$ is, in a sense, transported by itself.
\begin{lemma}\label{L:etaFlow}
	We have,
	\begin{align*}
		\prt_{t_1} \eta(t_1, t_2; \x) + \uu(t_1, \x) \cdot \grad \eta(t_1, t_2; \x)
			= 0.
	\end{align*}
\end{lemma}
\begin{proof}
From \cref{e:GroupComp},  $\y = \eta(t_1, t_2; \eta(t_2, t_1; \y))$; taking $d/d t_1$ of both sides of this identity,
\begin{align*}
	0
		&= \diff{}{t_1} \eta(t_1, t_2; \eta(t_2, t_1; \y)) \\
		&= \prt_1 \eta(t_1, t_2; \eta(t_2, t_1; \y))
			+ \grad \eta(t_1, t_2; \eta(t_2, t_1; \y)) \prt_{t_1} \eta(t_2, t_1; \y) \\
		&= \prt_1 \eta(t_1, t_2; \eta(t_2, t_1; \y))
			+ \uu(t_1, \eta(t_2, t_1; \y)) \cdot \grad \eta(t_1, t_2; \eta(t_2, t_1; \y)).
\end{align*}
Setting $\x = \eta(t_2, t_1; \y)$ gives the result.
\end{proof}

For any $(t, \x) \in \R \times \R^3$ let
\begin{itemize}
	\item
		$\pos(t, \x)$ be the point on $\Gamma_+$ at which the flow line
		through $(t, \x)$ intersects with $\Gamma_+$;
		
	\item
		let $\time(t, \x)$ be the time at which that intersection occurs.
\end{itemize}
For all $\x \in \Omega$, $\time(t, \x) \le t$. Because we extended $\uu$, $\time$ will always be defined, as long as we allow $\time(t, \x) = - \iny$, but $\pos(t, \x)$ is defined only when $\time(t, \x)$ is finite, and is meaningful only when $\time(t, \x) \ge 0$.

\TwoColumn{0.675}{0.3}{l}{c}
{
	We define the hypersurface,
	\begin{align*}
		S := \set{(t, \x) \in \ol{Q} \colon \time(t, \x) = 0}
	\end{align*}
	and the open sets $U_\pm \subset Q$,
	\begin{align*}
		U_-
			&:= \set{(t, \x) \in (0, T) \times \Omega
					\colon \time(t, \x) < 0}, \\
		U_+
			&:= \set{(t, \x) \in (0, T) \times \Omega
					\colon \time(t, \x) > 0}.
	\end{align*}
	For $t \in [0, T]$, we also define the sections, $S(t)$ and $U_\pm(t)$:
	\begin{align*}
		S(t) := \set{\x \colon (t, \x) \in S}, \quad
		U_\pm(t) := \set{\x \colon (t, \x) \in U_\pm}.
	\end{align*}
	This is illustrated in 2D in Figure 3.1.
}
{
\begin{center}
\scalebox{.34} % Change this value to rescale the drawing.
{
\hspace*{-135mm}
\begin{pspicture}(0,1.679176)(23.791176,13.383732)
\psbezier[linecolor=black, linewidth=0.04](13.918924,7.9412827)(12.748709,7.6313386)(12.587746,7.322413)(12.406923,6.460433081227662)(12.226101,5.598453)(12.621885,4.320358)(13.2028675,3.3135448)(13.783849,2.306732)(16.555172,1.7618151)(18.450235,1.7109543)(20.345299,1.6600934)(22.225084,2.6996412)(22.701765,4.9023185)(23.178444,7.1049953)(22.76421,10.290721)(21.767479,11.443267)(20.770746,12.595813)(18.353546,13.3280735)(17.205326,13.352278)(16.057104,13.376482)(14.091258,13.112)(13.999951,11.524082)(13.908645,9.936164)(14.356587,10.380652)(14.726,9.323452)(15.095413,8.266252)(14.488384,8.15939)(13.915904,7.941993)
\psbezier[linecolor=black, linewidth=0.04](16.086132,7.385118)(15.768206,7.200503)(14.984493,6.737444)(14.753847,6.306729312788421)(14.5232,5.8760147)(14.717874,4.9532375)(15.058672,4.613069)(15.399469,4.272901)(16.72482,4.1655045)(18.028118,4.159232)(19.331415,4.1529593)(20.433699,5.6526794)(20.54297,6.3248987)(20.65224,6.997118)(20.917751,9.356607)(20.067684,10.013475)(19.217617,10.670342)(17.886555,11.34138)(17.092356,11.298276)(16.298157,11.255172)(16.070095,10.907072)(16.107195,10.306666)(16.144295,9.70626)(16.936281,9.333769)(16.823126,8.708227)(16.709972,8.0826845)(16.41742,7.5504208)(16.076782,7.4057183)
\psellipse[linecolor=black, linewidth=0.04, dimen=outer](0.01,5.842278)(0.01,0.01)
\psrotate(18.523235, 7.3565426){-35.27809}{\psellipse[linecolor=black, linewidth=0.04, dimen=outer](18.523235,7.3565426)(0.735,1.0097059)}
\psline[linecolor=black, linewidth=0.08, arrowsize=0.05291667cm 2.0,arrowlength=1.4,arrowinset=0.0]{->}(15.105883,2.1821308)(15.317647,2.8174248)
\psline[linecolor=black, linewidth=0.08, arrowsize=0.05291667cm 2.0,arrowlength=1.4,arrowinset=0.0]{->}(13.305882,3.1586013)(13.8,3.5586014)
\psline[linecolor=black, linewidth=0.08, arrowsize=0.05291667cm 2.0,arrowlength=1.4,arrowinset=0.0]{->}(17.22353,1.8291895)(17.341177,2.4644837)
\psline[linecolor=black, linewidth=0.08, arrowsize=0.05291667cm 2.0,arrowlength=1.4,arrowinset=0.0]{->}(19.470589,1.7703661)(19.329412,2.452719)
\psline[linecolor=black, linewidth=0.08, arrowsize=0.05291667cm 2.0,arrowlength=1.4,arrowinset=0.0]{->}(14.235294,12.276248)(14.823529,11.946836)
\psline[linecolor=black, linewidth=0.08, arrowsize=0.05291667cm 2.0,arrowlength=1.4,arrowinset=0.0]{->}(21.37647,2.617425)(20.97647,3.135072)
\psline[linecolor=black, linewidth=0.08, arrowsize=0.05291667cm 2.0,arrowlength=1.4,arrowinset=0.0]{->}(22.551294,4.3148875)(21.919294,4.5670214)
\psline[linecolor=black, linewidth=0.08, arrowsize=0.05291667cm 2.0,arrowlength=1.4,arrowinset=0.0]{->}(22.905882,6.452719)(22.22353,6.499778)
\psline[linecolor=black, linewidth=0.08, arrowsize=0.05291667cm 2.0,arrowlength=1.4,arrowinset=0.0]{->}(22.8,8.60566)(22.117647,8.488013)
\psline[linecolor=black, linewidth=0.08, arrowsize=0.05291667cm 2.0,arrowlength=1.4,arrowinset=0.0]{->}(16.449339,13.324414)(16.538897,12.639847)
\psline[linecolor=black, linewidth=0.08, arrowsize=0.05291667cm 2.0,arrowlength=1.4,arrowinset=0.0]{->}(18.729412,13.076248)(18.517647,12.464483)
\psline[linecolor=black, linewidth=0.08, arrowsize=0.05291667cm 2.0,arrowlength=1.4,arrowinset=0.0]{->}(20.811764,12.217425)(20.411764,11.699778)
\psline[linecolor=black, linewidth=0.08, arrowsize=0.05291667cm 2.0,arrowlength=1.4,arrowinset=0.0]{->}(22.188236,10.746837)(21.576471,10.464483)
\psline[linecolor=black, linewidth=0.08, arrowsize=0.05291667cm 2.0,arrowlength=1.4,arrowinset=0.0]{->}(14.152941,10.264483)(14.741177,10.570366)
\psline[linecolor=black, linewidth=0.08, arrowsize=0.05291667cm 2.0,arrowlength=1.4,arrowinset=0.0]{->}(14.635294,8.299778)(15.082353,7.8056602)
\psline[linecolor=black, linewidth=0.08, arrowsize=0.05291667cm 2.0,arrowlength=1.4,arrowinset=0.0]{->}(12.694118,7.299778)(13.211764,6.9233074)
\psline[linecolor=black, linewidth=0.08, arrowsize=0.05291667cm 2.0,arrowlength=1.4,arrowinset=0.0]{->}(12.447059,5.0997777)(13.105883,5.288013)
\psline[linecolor=black, linewidth=0.1, arrowsize=0.05291667cm 2.0,arrowlength=1.4,arrowinset=0.0]{->}(19.329412,7.5586014)(18.811764,7.652719)
\psline[linecolor=black, linewidth=0.1, arrowsize=0.05291667cm 2.0,arrowlength=1.4,arrowinset=0.0]{->}(19.011765,6.770366)(18.588236,7.1468368)
\psline[linecolor=black, linewidth=0.1, arrowsize=0.05291667cm 2.0,arrowlength=1.4,arrowinset=0.0]{->}(17.870588,6.6409545)(18.317648,6.9468365)
\psline[linecolor=black, linewidth=0.1, arrowsize=0.05291667cm 2.0,arrowlength=1.4,arrowinset=0.0]{->}(18.470589,8.20566)(18.658823,7.7115426)
\psline[linecolor=black, linewidth=0.1, arrowsize=0.05291667cm 2.0,arrowlength=1.4,arrowinset=0.0]{->}(17.77647,7.5115423)(18.270588,7.323307)
\rput[bl](19.28,10.393895){\huge $S(t)$}
\rput[bl](18.17647,8.282131){\huge $\Gamma_{-}$}
\rput[bl](17.22353,4.9821306){\huge $U_{-}(t)$}
\rput[bl](19.752942,4.182131){\huge $U_{+}(t)$}
\rput[bl](22.441176,2.8644838){\huge $\Gamma_{+} = S(0)$}
\end{pspicture}
}

\small\textit{Figure 3.1}
\end{center}
}

	The hypersurface $S$ consists
	of all points $(t, \x)$ whose flow lines originated on $\Gamma_+$ at
	time zero; on $U_-$ these flow lines originated in $\Omega$ at time zero;
	on $U_+$ the flow lines originated on $\Gamma_+$ at positive time.
	Observe that $\pos(t, \x)$ is only meaningful on $\ol{U}_+$.

By virtue of \cref{e:GroupComp}, we have,
\begin{align}\label{e:mutDef}
	\begin{split}
		&\eta(\time(t, \x), t; \pos(t, \x)) = \x, \\
		&\eta(t, \time(t, \x); \x) = \pos(t, \x).
	\end{split}
\end{align}
Also note that $\time(t, \x) = t$ and $\pos(t, \x) = \x$ when $\x \in \Gamma_+$.

\begin{remark}\label{R:Brevity}
	We will often drop the $(t, \x)$ arguments on $\time$ and $\pos$
	for brevity.
\end{remark}

The regularity of $\eta$ given in \cref{L:etaRegularity} follows from entirely classical arguments,
so we omit its proof. Note that $\eta$ has one more derivative of time regularity (in both time variables) than $\uu$, making up for the loss of time regularity of $\uInputSpaceN$ from that of $\uSolSpaceN$.

\begin{lemma}\label{L:etaRegularity}
	The flow map
	$\eta \in C^{N + 1, \al}([0, T]^2 \times \R^3)$.
\end{lemma}

In \cref{L:Surface}, we show that the hypersurface $S$ has the regularity of the velocity field, and give conditions for when a function, regular on $U_\pm$ separately, can be glued together to obtain a regular function on all of $[0, T] \times \Omega$.

\begin{lemma}\label{L:Surface}
	The set $S$ is $C^{N + 1, \al}$ as a
	hypersurface in $[0, T] \times \R^3$.
	There is a $T^* > 0$ depending only upon $\norm{\uu}_{C^1(Q)}$ for which
	$\prt_t \time > 0$ on $\ol{U}_+$, while
	$\uu(t)$ remains transversal to $S(t)$ and
	$S(t) \subseteq \Omega$ for all $t \in [0, T^*]$.
	If $g \in C^{k, \al}(\ol{U}_-)$
	and $g \in C^{k, \al}(\ol{U}_+ \setminus \set{0} \times \Gamma_+)$
	for some $k \le N + 1$
	with $D^\beta g$ continuous on $S$ for all $\abs{\beta} = k$
	then $g \in C^{k, \al}(Q)$ with a norm no larger than
	the larger of its $C^{k, \al}$ norms on $U_\pm$.
\end{lemma}
\begin{proof}
	To simplify notation, we will make the argument as though $\Gamma_+$ has one component,
	the result for multiple components following simply by summing the estimates
	over each component.
	
	Because $\Gamma_+$ is a $C^{N + 1, \al}$ surface,
	we can write $\Gamma_+ := \varphi_0^{-1}(0)$ for some $\varphi_0 \in C^{N + 1, \al}(\R^3)$,
	which we can choose so that $\grad \varphi_0|_{\varphi_0^{-1}(0)} = \n$
	on $\Gamma_+$.
	Letting $\varphi_t(\x) = \varphi_0(\eta(t, 0; \x))$, we have
	$S(t) = \varphi_t^{-1}(0)$.
	Then by \cref{L:HolderComp},
	\begin{align*}
		\norm{\varphi_t}_{C^{N, \al}(\Omega)}
			&\le \norm{\varphi_0}_{C^{N, \al}(\R^3)}
				\brac{1 + \norm{\eta(t, 0; \x)}_{C^{N + 1}(\Omega)}}^{N + 1}.
	\end{align*}
	Since $S(t)$ is the surface $\Gamma_+$ transported to time $t$ by
	$\eta(0, \cdot; \cdot) \in C^{N + 1, \al}([0, T] \times \R^3)$,
	this shows that $S(t)$ is $C^{N + 1, \al}$ as a surface in $\Omega$.
	
	Now let $\varphi(t, \x) := \varphi_t(\x)$, so $\varphi \colon [0, T] \times \Omega \to \R$.
	Then $S = \varphi^{-1}(0)$ and we see, using also the regularity of $\eta$
	from \cref{L:etaRegularity}, that $S$ is $C^{N + 1, \al}$.
		
	Since $\varphi$ is transported by the flow, we have
	$
		\prt_t \varphi + \uu \cdot \grad \varphi = 0
	$. Since $\grad \varphi$ is normal to the surface $S(t)$, $\uu(t)$
	remains transversal to $S(t)$ as long as $\prt_t \varphi \ne 0$.
	But,
	\begin{align*}
		\abs{\prt_t \varphi|_{t = 0}}
			&= \abs{\uu_0 \cdot \grad \varphi_0|_{\varphi_0^{-1}(0)}}
			= \abs{\uu_0 \cdot \n}
			\ge U_{min}
			:= \min_{\Gamma_+} \abs{U^\n}
			> 0
	\end{align*}
	on $\Gamma_+$, and the regularity of $\uu$
	then ensures that
	$\abs{\prt_t \varphi(t)} > 0$, at least up to some finite time $T^* > 0$;
	hence, $\uu(t)$ remains transversal to $S(t)$
	and $\prt_t \time > 0$ on $\ol{U}_+$ for all $t \in [0, T^*]$. If necessary,
	we can always decrease $T^*$ so that
	$S(t) \subseteq \Omega$ for all $t \in [0, T^*]$.
					
	For the regularity of $g$, let $h := D^\beta g$ for any $\abs{\beta} = k$.
	Then $h \in C^\al(U_\pm)$ and is continuous on $S$. Hence, a simple application of
	the triangle inequality shows that we can glue $g|_{U_\pm}$ together
	along $S$ to obtain $h \in C^\al(Q)$, giving $g \in C^{k, \al}(Q)$.
	
	More explicitly,
	let $\y_\pm := (t_\pm, \x_\pm) \in \ol{U}_\pm$, and let $\y := (t, \x)$ be any point
	in $S$ along the line segment
	from $\y_-$ to $\y_+$.
	Let $h = D^\beta g$, which we note is defined everywhere
	on $S$. Then
	\begin{align*}
		&\frac{\abs{h(\y_-) - h(\y_+)}}{\abs{\y_- - \y_+}^\al}
			= \frac{\abs{h(\y_-) - h(\y) + h(\y)
				- h(\y_+)}}{\abs{\y_- -\y_+}^\al} \\
			&\qquad
			\le \frac{\abs{h(\y_-) - h(\y)}}{\abs{\y_- - \y_+}^\al}
					+ \frac{\abs{h(\y) - h(\y_+)}}{\abs{\y_- - \y_+}^\al} \\
			&\qquad
			= \la^\al \frac{\abs{h(\y_-) - h(\y)}}{\abs{\y_- -\y}^\al}
					+ (1 - \la)^\al
						\frac{\abs{h(\y) - h(\y_+)}}{\abs{\y - \y_+}^\al} \\
			&\qquad
			\le \la^\al \norm{h}_{\dot{C}^\al(U_-)}
				+ (1 - \la)^\al \norm{h}_{\dot{C}^\al(U_+)}
			\le \max \set{\norm{h}_{\dot{C}^\al(U_-)}, \norm{h}_{\dot{C}^\al(U_+)}},
	\end{align*}
	where
	\begin{align*}
		\la
			&= \frac{\abs{\y_- - \y}}{\abs{\y_- - \y_+}} \implies
		1 - \la  = \frac{\abs{\y_+ - \y}}{\abs{\y_- - \y_+}},
	\end{align*}
	where the equality for $1 - \la$ holds because we chose $\y$ to be colinear
	with $\y_\pm$.
	Similar estimates hold for all the norms.
\end{proof}

\begin{lemma}\label{L:muRegularity}
	Both $\time$ and $\pos$ are transported by the flow map for $\uu$; that is,
	\begin{align}\label{e:timeposTransported}
		\prt_t \time + \uu \cdot \grad \time = 0, \quad
		\prt_t \pos + \uu \cdot \grad \pos = 0,
	\end{align}
	and we have the identities
	(recall the definition of the operator $D$ in \cref{e:D}),
	\begin{align}\label{e:Dmu}
		\begin{split}
		\begin{array}{ll}
		\displaystyle
		\prt_t \time 
			= -\frac{\prt_{t_1} \eta(t, \time; \x) \cdot \n(\pos)}
				{U^\n(\time, \pos)},
		&\displaystyle
		\grad \time
			= -\frac{(\grad \eta(t, \time; \x))^T \n(\pos)}
				{U^\n(\time, \pos)}, \\
		\displaystyle
		D \time
			= -\frac{(D \eta(t, t_2; \x)|_{t_2 = \time})^T \n(\pos)}
				{U^\n(\time, \pos)}
		\end{array}
		\end{split}
	\end{align}
	and
	\begin{align}\label{e:Dpos}
		D \pos
			= D \eta(t, t_2; \x)|_{t_2 = \time} + \uu(\time, \pos) \otimes D \time.
	\end{align}
	Moreover, $\time$, $\pos$ lie in $C^{N + 1, \al}(\olUPlus)$.
\end{lemma}
\begin{proof}
	The quanitites $\time$ and $\pos$ are transported by the flow map for $\uu$ as
	in \cref{e:timeposTransported} because
	they are constant along flow lines by their definition.

	Taking the spatial derivatives, $\prt_{\x_\ell}$, of both sides of \cref{e:mutDef}$_2$ gives
	(recall \cref{R:Brevity})
	\begin{align*}
		\prt_{x_\ell} \eta^k(t, \time; \x)
			&= \prt_{t_2} \eta^k(t, \time; \x)  \prt_{x_\ell} \time
				+ \prt_{x_\ell} \eta^k(t, \time; \x) 
			= \prt_{x_\ell} \gamma^k,
	\end{align*}
	or, using \cref{e:FlowId}$_1$,
	\begin{align}\label{e:muRel1}
		\prt_{t_2} \eta(t, \time; \x) \otimes \grad \time
				+ \grad \eta(t, \time; \x)
			&= \uu(\time, \pos) \otimes \grad \time + \grad \eta(t, \time; \x)
			= \grad \pos.
	\end{align}
	
	Taking the time derivative of both sides of \cref{e:mutDef}$_2$ gives
	\begin{align}\label{e:muRel2}
		\prt_{t_1} \eta(t, \time; \x)
				+ \prt_t \time \prt_{t_2} \eta(t, \time; \x)
			= \prt_{t_1} \eta(t, \time; \x) + \prt_t \time \, \uu(\time, \pos)
			= \prt_t \pos.
	\end{align}
	
	We can write \cref{e:muRel1,e:muRel2} together in the form \cref{e:Dpos}.
	Now, $(D \pos)^T \n = 0$ because $\pos$ always lies on the boundary component $\Gamma_+$. 
	Thus (recall the comment following \cref{e:D}),
	\begin{align*}
		0
			&= (D \eta(t, t_2; \x)|_{t_2 = \time})^T \n
				+ (D \time \otimes \uu(\time, \pos)) \n
			= (D \eta(t, t_2; \x)|_{t_2 = \time})^T \n
				+ U^\n(\time, \pos)) D \time.
	\end{align*}
	(Note that $D \eta$ is a $4 \times 3$ matrix field.)
	From this, each of the expressions in \cref{e:Dmu} follow.

	Then \cref{e:Dmu} gives the regularity of $\time$,
	and \cref{e:Dpos} yields the regularity of $\pos$.
\end{proof}

\begin{lemma}\label{L:prt2gradeta}
	We have,
	\begin{align*}
		\prt_{t_2} \grad \eta(t_1, t_2; \z)
			&= \grad \uu(t_2, \eta(t_1, t_2; \z)) \grad \eta(t_1, t_2; \z)
				\text{ on } Q.
	\end{align*}
\end{lemma}
\begin{proof}
	Using \cref{e:FlowId},
	\begin{align*}
		\prt_{t_2} \grad \eta(t_1, t_2; \z)
			&= \prt_{t_2} \grad_\z \eta(t_1, t_2; \z)
			= \grad_\z \prt_{t_2} \eta(t_1, t_2; \z)
			= \grad_\z \uu(t_2, \eta(t_1, t_2; \z)) \\
			&= \grad \uu(t_2, \eta(t_1, t_2; \z)) \grad \eta(t_1, t_2; \z).
	\end{align*}
\end{proof}

%
% Section
%
\section{The pushforward}\label{S:Pushforward}

\noindent In this section we describe how to
extend the classical idea of the pushforward of the vorticity as a vector field to incorporate the generation of vorticity on the boundary. We start with a very brief overview of transport and the pushforward, specifically in flat space, but paying attention at the very beginning to issues of regularity. We then explain how we extend the pushforward to incorporate vorticity generation on the boundary, which we will use in \cref{S:Solutions} to produce a solution to \Our.

Our focus is on analysis, the regularity of our operations, rather than their geometric meaning. For a very readable treatment of what the pushforward means geometrically in the context of fluid mechanics we refer the reader to Sections 2.2 and 3.1 of \cite{BesseFrisch2017}.

\subsection*{Transport}
Before turning to the pushforward of a vector field, we review some basic facts regarding scalar transport---the pushforward of a scalar field. We define the transport operator applied to a scalar field $f(t, \x)$ by
\begin{align}\label{e:Lt}
	L_t f
		&:= \prt_t f + \grad_\uu f,
\end{align}
where $\grad_\uu$ is the directional derivative with respect to $\uu$. When $f$ has sufficient regularity, $\grad_\uu f = \uu \cdot \grad f$. As long as $f$ is, say, continuous, we can always write, for any fixed $s \in \R$,
\begin{align}\label{e:LtWeak}
	L_t (f(t, \eta(s, t; \x)))
		&= \diff{}{t} f(t, \eta(s, t; \x)).
\end{align}
We will also apply the transport operator to a vector field, component-by-component.

\subsection*{Pushforward}
The \textit{pushforward} of a vector field $\X_s$ by $\eta$
from time $s$ to time $t \in [0, T]$ is
\begin{align}\label{e:Y}
	(\eta(s, t)_* \X_s)(t, \x)
		:= \X_s(\eta(t, s; \x)) \cdot \grad \eta(s, t; \eta(t, s; \x)).
\end{align}
This presupposes that we stay within the domain of $\eta$ and $\X$.
We define the associated \textit{pushforward operator} $L$ by
\begin{align*}
	L \X := L_t \X - \X \cdot \grad \uu.
\end{align*}

\begin{lemma}\label{L:LY}
	For a vector field $\X \in C^{N, \al}(Q)$ for any $N \ge 0$, $L \X = 0$.
\end{lemma}
\begin{proof}
Holding $s$ and $\x$ fixed while applying \cref{e:LtWeak} to each component of $\X(t, x) := (\eta(s, t)_* \X_s)(t, \x)$ and appealing to \cref{L:prt2gradeta}, we have
\begin{align*}
	L_t (\X&(t, \eta(s, t; \x))
		= \diff{}{t} \X(t, \eta(s, t; \x)) \\
		&= \diff{}{t} \brac{\X_s(\eta(t, s; \eta(s, t; \x)))
			\cdot \grad \eta(s, t; \eta(t, s; \eta(s, t; \x)))} \\
		&= \diff{}{t} \brac{\X_s(\x)
			\cdot \grad \eta(s, t; \x)}
		= \X_s(\x) \cdot \diff{}{t} \grad \eta(s, t; \x) \\
		&= \X_s(\x) \cdot \pr{\grad \uu(t, \eta(s, t; \x)) \grad \eta(s, t; \x)} \\
		&= \pr{\X_s(\x) \cdot \grad \eta(s, t; \x)} \cdot \grad \uu(t, \eta(s, t; \x))
		= (\X \cdot \grad \uu)(t, \eta(s, t; \x)).
		\qedhere
\end{align*}	
\end{proof}
In geometric language, \cref{L:LY} tells us that the pushforward of a velocity field is Lie-transported (Lie-advected) by the flow (as in (3.9) of \cite{BesseFrisch2017}).

If $N \ge 1$, we can write the conclusion of \cref{L:LY} applied to a potential solution $\Y$ as
\begin{align*}
	\prt_t \Y + \uu \cdot \grad \Y = \Y \cdot \grad \uu,
\end{align*}
and we see the connection with \Our. The limitation, of course, is that this form of the pushforward does not account for inflow from the boundary (nor forcing, which we will consider later).

But if we let $\Y(t) = \eta(0, t)_* \Y_0$, we can see that the value of $\Y_0$ completely determines the value of $\Y$ on $U_-$. So other than the domain of $\Y$ shrinking in time as the flow sweeps $\Gamma_+$ through $\Omega$ (see Figure 3.1), the pushforward applies in an essentially classical way on $U_-$. Indeed, \cref{L:LY}, which is easily localized to $U_-$, gives that $L \Y = 0$ on $U_-$.

\subsection*{Inflow from the boundary}

By contrast, on $U_+$, we are assigned the value of $\Y = \H$ only on $\Gamma_+$, so that the values of $\Y$ on $U_+$ are entirely populated by the values of $\H$ on $[0, T] \times \Gamma_+$. Yet we wish the expression for the pushforward in \cref{e:Y} to continue to hold on $U_+$, as it yields \OurOne, 
but in the end we must connect the value of $\Y(t, \x)$ on $U_+$ to the time and place its value originated from on $\Gamma_+$. This is the main purpose of the functions $\time$ and $\pos$ that we defined in \cref{S:FlowMap}, and leads to the following definition:

\begin{definition}\label{D:Pushforward}
	Let $\X_0 \in C^\al(\Omega)$ and $\H \in C^\al([0, T] \times \Gamma_+)$.
	Define the pushforward of $\X_0$ by $\eta$ on $Q$
	with boundary value $\H$ by
	\begin{align*}
		\X(t, \x)
			&:=
			\begin{cases}
				(\eta(0, t)_* \X_0)(t, \x)
					&\text{on } U_-, \\
				(\eta(\time(t, \x), t)_* \H(\time(t, \x))(t, \x)
					&\text{on } U_+.
			\end{cases}
	\end{align*}
\end{definition}
Written out more fully, for $(t, \x) \in U_+$,
\begin{align*}
	\begin{split}
	\X(t, \x)
		&= \H(\time(t, \x), \pos(t, \x))
			\cdot \grad \eta(\time(t, \x), t; \eta(t, \time(t, \x); \x)) \\
		&= \H(\time(t, \x), \pos(t, \x))
			\cdot \grad \eta(\time(t, \x), t; \pos(t, \x)).
	\end{split}
\end{align*}

\begin{lemma}\label{L:YPlusPDE}
	With $\X$ as in \cref{D:Pushforward}, $L \X = 0$ on $U_+$.
\end{lemma}
\begin{proof}
	Let $(t, \x) \in U_+$ and fix $s \ge \time(t, \x)$.
	Because $\time$ and $\pos$ are constant along flow lines, so that
	$\time(t, \eta(s, t; \x)) = \time(s, \x)$ and
	$\pos(t, \eta(s, t; \x)) = \pos(s, \x)$, the verification of $L \X = 0$
	on $U_+$ parallels the calculation that led to \cref{L:LY}, though its
	expression is a little more complex.
	Holding $s$ and $\x$ fixed, we apply \cref{e:LtWeak} on $U_+$, to give
	\begin{align*}
		L_t (\X&(t, \eta(s, t; \x))
			= \diff{}{t} \X(t, \eta(s, t; \x)) \\
			&= \diff{}{t}
				\brac{
					\H(\time(t, \eta(s, t; \x)), \pos(t, \eta(s, t; \x)))
					\cdot \grad \eta(\time(t, \eta(s, t; \x)), t; \pos(t, \eta(s, t; \x)))
				}  \\
			&= \diff{}{t}
				\brac{
					\H(\time(s, \x), \pos(s, \x))
					\cdot \grad \eta(\time(s, \x), t; \pos(s, \x))
				} \\
			&= \H(\time(s, \x), \pos(s, \x))
				\cdot \prt_{t_2} \grad \eta(\time(s, \x), t; \pos(s, \x)) \\
			&= \H(\time(s, \x), \pos(s, \x)) \cdot
				\brac{\grad \uu(t, \eta(\time(s, \x), t; \pos(s, \x)))
					\grad \eta(\time(s, \x), t; \pos(s, \x))} \\
			&= \brac{\H(\time(s, \x), \pos(s, \x)) \cdot \grad \eta(\time(s, \x), t; \pos(s, \x))}
				\grad \uu(t, \eta(\time(s, \x), t; \pos(s, \x))) \\
			&= \brac{\H(\time(t, \eta(s, t; \x)), \pos(t, \eta(s, t; \x)))
					\cdot \grad \eta(\time(t, \eta(s, t; \x)), t; \pos(t, \eta(s, t; \x)))} \\
			&\qquad
				\grad \uu(t, \eta(\time(t, \eta(s, t; \x)), t; \pos(t, \eta(s, t; \x)))) \\
			&= \brac{\H(\time(t, \z), \pos(t, \z))
					\cdot \grad \eta(\time(t, \z), t; \pos(t, \z))}
				\grad \uu(t, \eta(\time(t, \z), t; \pos(t, \z)))|_{\z = \eta(s, t; \x)} \\
			&= (\X \cdot \grad \uu)(t, \eta(s, t; \x)).
	\end{align*}
	This shows that $L \X = 0$ on $U_+$.
	(Note that assuming higher regularity for $\H$ was not required for this proof.)
\end{proof}

%
% Section
%
\section{Lagrangian and Eulerian solutions}\label{S:Solutions}

\noindent To define a Lagrangian solution, we first explain how to handle forcing. We then show, in \cref{P:RegNFromS}, how to obtain the regularity of Lagrangian solutions from the compatibility condition $\cond_N$. Finally, we relate our Lagrangian solutions to Eulerian solutions in \cref{P:LinearExistenceFirstPart}.

\begin{assumption}\label{A:TStarLimit}
	$T^* = T$, where $T^*$ is as in \cref{L:Surface}.
\end{assumption}

\begin{remark}\label{R:TStarLimit}
	In \cref{R:BeyondTStar}, we show how to drop \cref{A:TStarLimit}.
\end{remark}

\subsection*{Forcing}

To treat forcing, which we have so far not considered, we use a version of Duhamel's principle.
We define $\G$ is as in (3.21) Chapter 4 of \cite{AKM},
\begin{align}
	\G(t, \x)
		&:= \int_{\ol{\time}(t, \x)}^t (\eta(s, t)_* \g(s))(t, \x) \, ds,
				\label{e:YForcing}
\end{align}
where $\ol{\time}(t, \x) = \max \set{0, \time(t, \x)}$.

\begin{prop}\label{P:Forcing}
	Assume that $\g \in C^\al(Q)$ and $\uu \in \uInputSpace$.
	Then $\G \in C^\al(Q)$ and $L \G = \g$ on $Q$ weakly.
	If $\g \in \mathring{C}^{N, \al}(U_\pm) \cap C^\al(Q)$ and $\uu \in \uInputSpaceN$ for $N \ge 1$ then 
	$\G \in C^{N, \al}(U_\pm)$ with
	\begin{align*}
		\prt_t \G(t, \x)
			&= \g(t, \x) - (\eta(\ol{\time}(t, \x), t)_* \g(\ol{\time}(t, \x)))(t, \x)
				\prt_t \ol{\time}(t, \x)
			+ \int_{\ol{\time}(t, \x)}^t
				\prt_t ((\eta(s, t)_* \g(s))(t, \x)) \, ds, \\
		\grad \G(t, \x)
			&=  - (\eta(\ol{\time}(t, \x), t)_* \g(\ol{\time}(t, \x)))(t, \x)
				\otimes \grad \ol{\time}(t, \x)
			+ \int_{\ol{\time}(t, \x)}^t
				\grad_x ((\eta(s, t)_* \g(s))(t, \x)) \, ds
	\end{align*}
	for all $(t, \x) \in Q$,
	noting that the terms involving $\prt_t \ol{\time}(t, \x)$ and $\grad \ol{\time}(t, \x)$
	have a potential singularity along $S$.
\end{prop}
\begin{proof}
	The expressions for $\prt_t \G$ and $\grad \G$ follow from
	applying the chain rule for integrals to \cref{e:YForcing},
	as does
	\begin{align*}
		L \G
			&= (\eta(t, t)_* \g(t))(t, \x)
				- (\eta(\ol{\time}(t, \x), t)_* \g (\ol{\time}(t, \x)))(t, \x)
					\pr{\prt_t + \uu \cdot \grad} \ol{\time}(t, \x)  \\
			&\qquad
				+ \int_{\ol{\time}(t, \x)}^t
					L (\eta(s, t)_* \g(s))(t, \x) \, ds \\
			&= \g(t, \x)
				- (\eta(\ol{\time}(t, \x), t)_* \g (\ol{\time}(t, \x)))(t, \x)
					\pr{\prt_t + \uu \cdot \grad} \ol{\time}(t, \x),
	\end{align*}
	since $L (\eta(s, t)_* \g(s))(t, \x) \equiv 0$ for all $s$.
	Now, $\ol{\time}$ either equals 0 or $\time(t, \x)$, but either way, by
	virtue of \cref{L:muRegularity}, we see
	that its value is transported along flow lines. Hence, also
	$\pr{\prt_t + \uu \cdot \grad} \ol{\time}(t, \x) = 0$, and we conclude that
	$L \G = \g$---though only weakly because  $\prt_t \G$ and $\grad \G$
	are discontinuous along $S$.
	The regularity of $\G$ follows in a very classical way, employing the lemmas
	in \cref{S:HolderLemmas}.
\end{proof}

Define,
\begin{align}
	\pos_0
		&= \pos_0(t, \x)
		:= \eta(t, 0; \x)
		\text{ on } \ol{U}_-,
			\label{e:pos0} \\
	B_-
		&= B_-(t, \x)
		:= \grad \eta(0, t; \eta(t, 0; \x))
		= \grad \eta(0, t; \pos_0)
		\text{ on } \ol{U}_-, \label{e:Bminus} \\
	B_+
		&= B_+(t, \x)
		:= \grad \eta(\time, t; \eta(t, \time; \x))
		= \grad \eta(\time, t; \pos)
		\onUPlus. \label{e:Bplus}
\end{align}

\begin{definition}[Lagrangian solution to \Our]\label{D:LagrangianSolution}
Define $\Y$ by $\Y|_{U_\pm} = \Y_\pm$, where
\begin{align}\label{e:LagrangianForm}
	\begin{split}
	\Y_-(t, \x)
		&:= B_- \Y_0(\pos_0) + \G_-(t, \x), \\
	\Y_+(t, \x)
		&:= B_+ \H(\time, \pos) + \G_+(t, \x), \\
	\G_-(t, \x)
		&:= \int_0^t (\eta(s, t)_* \g(s))(t, \x) \, ds, \\
	\G_+(t, \x)
		&:= \int_{\time(t, \x)}^t (\eta(s, t)_* \g(s))(t, \x) \, ds.
	\end{split}
\end{align}
We say that $\Y$ is the Lagrangian solution to \Our.
\end{definition}

\begin{lemma}\label{L:YRegularityOnUpm}
	Assuming data regularity $N \ge 0$, $\Y_\pm \in C^{N, \al}(U_\pm)$.
\end{lemma}
\begin{proof}
	This follows directly from the definition
	of the pushforward, the regularity assumed on $\Y_0$ and $\H$,
	the regularity of $\eta$, $\time$, $\pos$, and $\G$
	given by \cref{L:etaRegularity,L:muRegularity,P:Forcing}, and
	the lemmas in \cref{S:HolderLemmas}.
\end{proof}

We will obtain uniqueness of solutions later in \cref{P:LinearExistenceFirstPart}, but we will find a need in the proof of \cref{P:RegNFromS} for a weaker kind of uniqueness (that is, with stronger assumptions) that allows us to know that certain Eulerian solutions are, in fact, Lagrangian solutions. That is the purpose of \cref{L:OnUPm}.

\begin{lemma}\label{L:OnUPm}
	Assume that the data has regularity $N \ge 1$, $\Y \in C^{1, \al}(U_\pm)$
	and $\Y \in C^\al(Q)$, $\Y = \H$ on $[0, T] \times \Gamma_+$, $\Y(0) = \Y_0$ on $\Omega$
	for some $\H \in C^{N, \al}([0, T] \times \Gamma_+)$, and $\Y_0 \in \vortInitSpaceN$.
	If
	\begin{align}\label{e:YWant}
 		\prt_t \Y + \uu \cdot \grad \Y - \Y \cdot \grad \uu
			&= \g
			\text{ on } U_\pm
	\end{align}
	then $\Y$ is the Lagrangian solution given by \cref{e:LagrangianForm}.
\end{lemma}
\begin{proof}
	From \cref{L:LY,L:YPlusPDE,L:YRegularityOnUpm} the Lagrangian solution
	given by \cref{e:LagrangianForm} satisfies \cref{e:YWant}, so we need
	only show uniqueness of solutions to \cref{e:YWant} with
	the same values of $\H$ and $\Y(0)$.
	Since \cref{e:YWant} is linear,
	this is equivalent to showing that given vanishing values of $\H$, $\Y(0)$,
	and forcing, the only solution is $\Y \equiv 0$.
	
	We suppose, then, that
	\begin{align}\label{e:OurLinearZero}
		\begin{cases}
			\prt_t \Y + \uu \cdot \grad \Y - \Y \cdot \grad \uu = 0
				&\text{in } U_+ \cup U_-, \\
			\Y = 0
				&\text{on } [0, T] \times \Gamma_+, \\
			\Y(0) = 0
				&\text{on } \Omega.
		\end{cases}
	\end{align}

	Multiplying by $\Y$ and integrating over $U_+ \cup U_-$, we have
	\begin{align*}
		\diff{}{t} &\norm{\Y}_{L^2(\Omega)}^2
			= - 2 (\uu \cdot \grad \Y, \Y)_{L^2(U_+(t))}
				- 2 (\uu \cdot \grad \Y, \Y)_{L^2(U_-(t))}
				+ 2 (\Y \cdot \grad \uu, \Y)_{L^2(\Omega)},
	\end{align*}
	using that $S(t)$ has measure zero.
	But, 
	\begin{align*}
		- 2 (\uu \cdot \grad \Y, &\Y)_{L^2(U_+(t))}
				- 2 (\uu \cdot \grad \Y, \Y)_{L^2(U_-(t))}
				\\
		&= - (\uu, \grad \abs{\Y}^2)_{L^2(U_+(t))} - (\uu, \grad \abs{\Y}^2)_{L^2(U_-(t))} \\
			&= - \int_{S(t)} (\uu \cdot \n) \abs{\Y}^2 + \int_{S(t)} (\uu \cdot \n) \abs{\Y}^2
				- \int_{\Gamma_-} U^\n \abs{\Y}^2
			\le 0.
	\end{align*}
	We used here that $\Y \in C^\al(Q)$, so the two boundary integrals over $S(t)$,
	the shared portion of $\prt U_+(t)$ and $\prt U_-(t)$, properly oriented cancel.
	Moreover, we used that
	$\Y$ has sufficient regularity to integrate the term
	$(\uu, \grad \abs{\Y}^2)$ by parts over $U_+$ and $U_-$ separately.
	Using that $(\Y \cdot \grad \uu, \Y)_{L^2(\Omega)} \le
	\norm{\grad \uu}_{L^\iny(\Omega)} \norm{\Y}_{L^2(\Omega)}^2$ and
	applying \Gronwalls lemma yields $\Y \equiv 0$ in $Q$,
	establishing uniqueness.
\end{proof}

\subsection*{Regularity across $S$}

The solution can be decomposed as in \cref{e:LagrangianForm} because \Our is a linear problem, the velocity field $\uu$ and so the flow map $\eta$ being given.
We see that formally \Our holds; however, this applies only distributionally on $U_\pm$ separately until we can establish sufficient regularity of $\Y$ across $S$, as we will see in the proof of \cref{P:RegNFromS}. Before turning to its proof, however, let us consider some of the difficulties in treating it.

For transport or transport-stretching equations, much of the analysis is done in more-or-less Lagrangian form, following the solution along flow lines. There, the natural operator that emerges is $L_t$ or $L$. The regularity of $L_t \Y_0$ is obtained quite easily in this way, but it does not allow the regularity of $\prt_t \Y$ and $\uu \cdot \grad \Y$ to be treated separately. Indeed, even without inflow, $\Y_0$ discontinuous along, say, a curve is no obstacle to obtaining solutions in which $L \Y_0 = 0$, so we should not expect to be able to separate them.

Nonetheless, classically, the regularity of the pushforward of $\prt_t \Y$ and so, as we shall see, $\grad \Y$ is fairly easily tied directly to the regularity of $\Y_0$. The difficulty we face, is that while $\Y_-$ is a classical pushforward off a domain from time zero, $\Y_+$ is a pushforward off of a surface into a domain, and the behavior of higher time derivatives of $\Y_+$ on the interface $S$ is not as easily connected to its behavior at time zero or to the behavior of $\Y_\pm$ on $S$.

Closely related to this is that for $N \ge 1$ and nonzero forcing, solutions $\Y \in \vortSpaceN$ do not require that both terms in the decomposition of $\Y_\pm$ in \cref{e:LagrangianForm} be in $\vortSpaceN$. Indeed, we can see from \cref{P:Forcing} that $\G$ will never have higher than $\vortSpace$ regularity unless we impose the strong condition that $\prt_t^j \g(0)$ vanish on $\Gamma_+$ for all $1 \le j \le N$.

In the proof of \cref{P:RegNFromS} we reference \cref{P:N1Condition}, which we prove in the next section, because it is of a different character than the rest of the proof.

\begin{prop}\label{P:RegNFromS}
	Assume that the data has regularity $N$ as in \cref{D:NReg} for some $N \ge 0$
	and let $\Y$ be the Lagrangian solution to \Our as in \cref{D:LagrangianSolution}.
	Then $\Y \in \vortSpaceN$ if and only if $\cond_N$ holds.
\end{prop}
\begin{proof}
	First observe that $\Y_\pm \in C^{N, \al}(U_\pm)$ follows from \cref{L:YRegularityOnUpm}
	without the need for any additional conditions.
	Also, $\Y \in \vortSpaceN$ gives cond$_N$,
	as we can most readily see from the form of $\cond_N$ in \cref{e:CondNAlt}.
	
	We assume, then, that $\cond_N$ holds and will prove that
	$\Y \in \vortSpaceN$.
		
	Assume first that $N = 0$. By \cref{P:Forcing}, $G \in C^\al(Q)$, and we can see
	from \cref{e:LagrangianForm} that
	along the hypersurface $S$, where $\time(t, \x) = 0$ and $B_+ = B_-$, the two expressions
	for $\Y(t, \x)$ agree if and only if $\H(0) = \Y_0$ on $\Gamma_+$.
	It follows from \cref{L:Surface} that $\Y \in C^\al(Q)$,
	completing the proof for $N = 0$.
		
	Now assume that $N = 1$.
	As already observed, $\Y \in C^{1, \al}(U_\pm)$. By
	\cref{L:Surface}, then, we see that $\prt_t \Y \in C^\al(Q)$
	if and only if $\prt_t \Y_- = \prt_t \Y_+$ along $S$
	and $\grad \Y \in C^\al(Q)$
	if and only if $\grad \Y_- = \grad \Y_+$ along $S$:
	both conditions hold by \cref{P:N1Condition} since we assumed that $\cond_1$ holds.
	Hence, $\Y \in C^{1, \al}(Q)$, and we can write
	$\prt_t \Y + \uu \cdot \grad \Y = \Y \cdot \grad \uu + \g$ on all of $Q$.
	This completes the proof for $N = 1$.
	
	Now assume $N = 2$ data regularity. Then $\Y \in C^{2, \al}(U_\pm)$,
	as already observed for general $N$, and by the $N = 1$ result we know that
	$\Y \in C^{1, \al}(Q)$.
		
	Next, let $\ZZ := \prt_t \Y$, so $\ZZ \in C^{1, \al}(U_\pm)$ and
	$\ZZ \in C^\al(Q)$. Then
	\begin{align}\label{e:ZEq}
 		\prt_t \ZZ + \uu \cdot \grad \ZZ - \ZZ \cdot \grad \uu
			&= \h
			:= \prt_t \g - \prt_t \uu \cdot \grad \Y + \Y \cdot \grad \prt_t \uu
				\in \mathring{C}^{1, \al}(U_\pm) \cap C^\al(Q).
	\end{align}
	Equality holds classically on $U_\pm$, so, in fact, $\ZZ$ is the Lagrangian solution
	as in \cref{D:LagrangianSolution} with forcing $\h$,
	initial value $\prt_t \Y(0)$,
	and $\ZZ = \prt_t \H$ on $[0, T] \times \Gamma_+$.
	This follows from \cref{L:OnUPm} or,
	with much greater trouble, directly from applying $\prt_t$ to the expression
	for $\Y$ in \cref{e:LagrangianForm}.

	We can see that $\cond_2$ for $\Y$ with forcing $\g$ and boundary value $\H$
	is $\cond_1$ for $\ZZ$ with forcing $\h$ and boundary value $\prt_t \H$,
	so \cref{P:N1Condition} gives that
	$
		\prt_t \ZZ_- = \prt_t \ZZ_+ \onS
	$
	and
	$
		\grad \ZZ_- = \grad \ZZ_+ \onS
	$.
	Applying \cref{L:Surface}, this, in turn, gives $\prt_t \Y = \ZZ \in C^{1, \al}(Q)$.
		
	Now, we cannot immediately apply this same argument to $\W := \prt_j \Y$,
	since $\cond_2$ for $\Y$
	does not become $\cond_1$ for $\W$. But we can take $\prt_j$ of \OurOne, giving
	\begin{align}\label{e:Weq}
 		\prt_t \W + \uu \cdot \grad \W - \W \cdot \grad \uu
			&= \jj
			:= \prt_j \g - \prt_j \uu \cdot \grad \Y + \Y \cdot \grad \prt_j \uu
			\in \mathring{C}^{1, \al}(U_\pm) \cap C^\al(Q).
	\end{align}
	Since we already know that $\prt_t \W = \prt_t \prt_j \Y = \prt_j \ZZ \in C^\al(Q)$, we have
	$\uu \cdot \grad \W \in C^\al(Q)$. But $\cond_0$ and the
	$N = 2$ data regularity give that $\grad \W$ is $C^\al$-continuous
	along $S(t)$. Letting $Q^* = (0, T^*) \times \Omega$,
	where $T^* > 0$ is given by \cref{L:Surface},
	the transversality of $\uu \in C^{3, \al}(Q^*)$ across $S(t)$ gives that
	$\grad \W$ is $C^\al$-continuous in the direction perpendicular
	to $S(t)$. We conclude that  
	$\grad \W = \grad \prt_j \Y \in C^\al(Q^*)$
	and so $\grad^2 \Y \in C^\al(Q^*)$.
	Combined with the regularity of $\ZZ$, we have $\Y \in C^{2, \al}(Q^*)$

	We now have the regularity to return to \cref{e:Weq} and view it as the solution to
	\begin{align}\label{e:OurWLinear}
		\begin{cases}
			\prt_t \W + \uu \cdot \grad \W - \W \cdot \grad \uu = \jj
				&\text{in } Q, \\
			\W = \J
				&\text{on } [0, T] \times \Gamma_+, \\
			\W(0) = \prt_j \Y_0
				&\text{on } \Omega,
		\end{cases}
	\end{align}
	where $\J := \prt_j \Y$. Since we know that $\W \in C^{1, \al}(Q^*)$,
	it follows from \cref{P:N1Condition} that $\cond_1$ holds for $\W$.
	But it then follows again from \cref{P:N1Condition} that, in fact,
	$\W \in C^{1, \al}(Q)$. Combined with the regularity of $\ZZ$,
	we can conclude that $\Y \in C^{2, \al}(Q)$. 

	This argument inducts to any value of $N \ge 3$.
\end{proof}

In the proof of \cref{P:RegNFromS}, we reduced the $N \ge 2$ case inductively to the $N = 1$ case, using the convenient form of $\cond_N$ in \cref{e:CondNAlt}. Formally, we could reduce all the way to the $N = 0$ case (removing the need for \cref{S:N1RegCond} altogether). In \cref{e:ZEq}, however, $\ZZ$ would only lie in $C^\al(U_\pm)$ and $\h$ at best would lie in the negative \Holder space $C^{\al - 1}(Q)$, which we do not have the tools to handle.

It is not hard to see that $\cond_0$ is the Rankine-Hugoniot condition that allows $\Y$ as given by \cref{D:LagrangianSolution} to be a weak Eulerian solution to \Our. There is a somewhat delicate issue regarding regularity, however, as in the classical case one typically considers $C^1$ solutions, while for $N = 0$, $\Y$ is only $C^\al$. Therefore, we present a proof in \cref{L:RankineHuginoit}, which we will need to  obtain Eulerian solutions in \cref{P:LinearExistenceFirstPart}.

\begin{lemma}\label{L:RankineHuginoit}
	Assume data regularity $N = 0$ and let $\Y$ be the Lagrangian solution to \Our as in
	\cref{D:LagrangianSolution}. If $\cond_0$ holds then $\Y$ is a weak Eulerian solution
	to \Our.
\end{lemma}
\begin{proof}
Let $(\uu_i)$ be the sequence in $\uBoundarySpaceOne{U_i^\n}$ approximating $\uu \in \uInputSpace$ given by \cref{L:VelDense01}. Also let $(\Y_{i, 0})$ be a sequence in $\vortInitSpaceOne$ with $\Y_{i, 0} \to \Y_0$ in $\vortInitSpace$. Let $\Y_i$ be the Lagrangian solution to \Our as in \cref{D:LagrangianSolution} with $\uu_i$ in place of $\uu$ and $\Y_{i, 0}$ in place of $\Y_0$ (leaving $\H$ and $\g$ unchanged). In what follows, $i$ subscripts on quantities, such as $U_{i, \pm}$, will refer to the corresponding quantity for $\uu$, $\Y$ applied with $\uu_i$, $\Y_i$.

Neither $\cond_0$ nor $\cond_1$ need hold for $\Y_i$, so we only have $\Y_{i, \pm} \in C^{1, \al}(U\pm)$, but \cref{L:LY,L:YPlusPDE,L:YRegularityOnUpm} do give, as in the proof of \cref{{L:OnUPm}}, that
\begin{align*}
	\prt_t \Y_i + \uu_i \cdot \grad \Y_i - \Y_i \cdot \grad \uu_i = \g
		\text{ on } U_{i, \pm}.
\end{align*} 
Let $\varphi \in \Cal{D}(Q)$. Multiplying the above equation by $\varphi$ and integrating by parts gives
\begin{align*}
	((\Y_i, \prt_t \varphi))
			+ (&(\Y_i \cdot \grad \uu_i, \varphi))
			+ ((\g, \varphi))
			+ \int_0^T \int_{S_i(t)} (\uu_i \cdot \n)
				(\Y_{i, -} - \Y_{i, +}) \cdot \varphi
		= 0,
\end{align*}
where $((\cdot, \cdot))$ is the pairing in $\Cal{D}'(Q), \Cal{D}(Q)$ and where we have arbitrarily chosen $\n$ on $S_i(t)$
to be outward from $U_{i, +}(t)$.

Because $\G_i \in C^\al(Q)$ by \cref{P:Forcing}, we have $\abs{\G_{i, +} - \G_{i, -}} \to 0$ on $S_i(t)$. Also,
\begin{align*}
	B_{i, -} \Y_{i, 0}(\pos_{i, 0}) - B_{i, +} \H(\time_i, \pos_i)
		&= B_{i, -} (\Y_{i, 0}(\pos_{i, 0}) - \H(\time_i, \pos_i))
			+ (B_{i, -} - B_{i, +}) \H(\time_i, \pos_i) \\
		&= B_{i, -} (\Y_{i, 0}(\pos_i) - \H(0, \pos_i))
			\text{ on } S_i(t),
\end{align*}
since on $S_i(t)$, $\time_i = 0$, $\pos_{i, 0} = \pos_i$, and $B_{i, -} = B_{i, +}$, without the need for $\cond_0$ to hold. But $\Y_{i, 0} \to \Y_0$ in $C^\al(\Omega)$ and $\Y_0 = \H$ on $\Gamma_+$, and $\pos_i(t, \x) \in \Gamma_+$ for $(t, \x) \in S$, so we see that
\begin{align*}
	B_{i, -} \Y_{i, 0}(\pos_{i, 0}) - B_{i, +} \H(\time_i, \pos_i) \to 0
		\text{ on } S_i(t).
\end{align*}
Hence,
$
	\Y_{i, -} - \Y_{i, +}
		\to 0
		\text{ on } S_i(t).	
$
We conclude from $\uu_i \to \uu$ in $\uInputSpace$ that in the limit,
\begin{align*}
	((\Y, \prt_t \varphi))
			+ (&(\Y \cdot \grad \uu, \varphi))
			+ ((\g, \varphi))
		= 0,
\end{align*}
which shows that $\Y$ is a weak Eulerian solution to \Our.
\end{proof}

\begin{prop}\label{P:LinearExistenceFirstPart}
	Assume that the data has regularity $N$ as in \cref{D:NReg} for some $N \ge 0$
	and let $\Y$ be the Lagrangian solution to \Our as in \cref{D:LagrangianSolution}.
	Assume that $\cond_N$ holds.
	The estimate in \cref{e:RegBound} holds
	and $\Y$ is a weak Eulerian solution as in \cref{D:Weak}.
	Moreover, $\Y$ is the unique classical solution in $C^{N, \al}(Q)$ to \Our when
	$N \ge 1$.
\end{prop}
\begin{proof}
	Since $\cond_N$ holds, we know from \cref{P:RegNFromS} that $\Y$ is
	a Lagrangian solution to \Our.
	The bound in \cref{e:RegBound} follows from the Lagrangian
	form of $\Y$ in \cref{e:LagrangianForm}, applying \cref{L:Surface} and
	\cref{P:Forcing}. That $\Y$ is a weak Eulerian solution for $N = 0$
	follows from \cref{L:RankineHuginoit}.
	For $N \ge 1$, all the terms in \Our are defined classically, so $\Y$ is a
	classical and so weak solution to $\Our$.
	
	To prove uniqueness of Eulerian solutions for $N \ge 1$, suppose that $\Y_1$ and $\Y_2$
	are two solutions to \Our. Letting $\Y = \Y_1 - \Y_2$, we see that
	\begin{align*}
		\begin{cases}
			\prt_t \Y + \uu \cdot \grad \Y - \Y \cdot \grad \uu = 0
				&\text{in } Q, \\
			\Y = 0
				&\text{on } [0, T] \times \Gamma_+, \\
			\Y(0) = 0
				&\text{on } \Omega.
	\end{cases}
	\end{align*}
	Since $N \ge 1$, we have enough regularity to
	multiply by $\Y$, and integrate over $\Omega$ using that
	$\Y = 0$ on the inflow boundary, to obtain
	\begin{align}\label{e:YUnique}
		\begin{split}
		\frac{1}{2} \diff{}{t} \norm{\Y}_{L^2}^2
			&\le \norm{\grad \uu}_{L^\iny(Q)} \norm{\Y}_{L^2}^2
				- (\uu \cdot \grad \Y, \Y)
			\le C \norm{\Y}_{L^2}^2 - \frac{1}{2} (\uu, \grad \abs{\Y}^2) \\
			&= C \norm{\Y}_{L^2}^2 
				- \frac{1}{2} \int_{\Gamma_+} U^\n \abs{\Y}^2
				- \frac{1}{2} \int_{\Gamma_-} U^\n \abs{\Y}^2
			\le C \norm{\Y}_{L^2}^2,
		\end{split}
	\end{align}
	since $\abs{\Y}^2 = 0$ on $\Gamma_+$ and $U^\n > 0$ on $\Gamma_-$.
	Applying \Gronwalls lemma gives $\Y = 0$ so $\Y_1 = \Y_2$.
\end{proof}

\begin{remark}\label{R:BeyondTStar}
	We have obtained the results of this section
	under \cref{A:TStarLimit}
	or, equivalently, we have obtained these results
	up to time $T^*$.
	To extend beyond $T^*$, we reset time so that
	$T^*$ becomes time zero.
	Because compatibility conditions hold for solutions
	at all positive times, we know that $\cond_N$ holds at the new time zero,
	and extend the result for another $T^*$ time interval, repeating this
	process as needed.
	This allows us to drop \cref{A:TStarLimit}.
\end{remark}

%
% Section
%
\section{Regularity condition for $N = 1$}\label{S:N1RegCond}

\noindent In this section we prove \cref{P:N1Condition}, which is the key to obtaining $\Y \in \vortSpaceN$.

\begin{prop}\label{P:N1Condition}
	Assume the data has $N = 1$ regularity, though we only assume that
	$\g \in \mathring{C}^{1, \al}(U_\pm) \cap C^\al(Q)$, and assume that $\cond_0$ holds.
	Let $\Y$ be the Lagrangian solution as in \cref{D:LagrangianSolution}.
	Then
	$
		D \Y_- = D \Y_+ \onS
			\text{ if and only if } \cond_1
	$
	holds.
	(Recall that $D$ is defined in \cref{e:D}.)
\end{prop}

The direct calculation showing that $D \Y_+ = D \Y_-$ on $S$ if $\cond_1$ holds is quite lengthy. Instead, we use the indirect method below. But first, we develop some necessary tools.

We know from the proof of \cref{P:RegNFromS} that regularity of $\Y$ on $Q$ derives from its regularity across the hypersurface $S$. To assess such regularity, we would like to be able to treat $\Y_\pm$ uniformly on the same domain, but their domains overlap only on $S$. Therefore, we extend $\Y_-$ to all of $Q$ by extending $\Y_0$ to $C^{1, \al}(\R^3)$ and $\g$ to $C^\al([0, T] \times \R^3)$. Then because, as in \cref{S:FlowMap}, we have extended $\uu$ and so $\eta$, the expression for $\Y_-$ in \cref{e:LagrangianForm} is defined on all of $Q$ and lies in $C^{1, \al}(Q)$ as long as $\uu \in \uInputSpaceOne$.
Moreover, it follows that
$
	\prt_t \Y_- + \uu \cdot \grad \Y_- - \Y_- \cdot \grad \uu = \g
$
on $Q$.

Because $S$ is at least a $C^{2, \al}$ hypersurface in $Q$ by \cref{L:Surface}, one-sided derivatives of $\Y$ in time and space up to order $N = 1$ exist on $U_\pm$. Hence, $\Y_- \in C^{1, \al}(\ol{U}_-)$ with $\Y_+ \in C^{1, \al}(\olUPlus)$. We remove $\set{0} \times \Gamma_+$ from $\ol{U}_+$, since there are no spatial derivatives of $\Y_+$ normal to the boundary at time zero. The time derivatives, however, are defined and continuous on all of $\ol{U}_+$.

In any case, one-sided derivatives of $\Y_\pm$ exist on $S$, so we can freely take derivatives of both $\Y_-$ and $\Y_+$ on $\olUPlus$. Restricted to $S$, any calculation for $\Y_-$ that depends only upon its $C^{1, \al}(\ol{U}_-)$ regularity is independent of the extension when restricted to $S$. Hence, the calculations on $S$ that follow do not depend upon the manner in which we extend $\Y_0$ and $\g$. These extensions also allow us to make sense of $\pos_0$, $B_-$, and $B_+$ of \cref{e:pos0,e:Bminus,e:Bplus} on $\olUPlus$, along with the matrix-valued function,
\begin{align}\label{e:M}
	M &= M(t, \x) := B_-^{-1} - B_+^{-1}
		\onUPlus. 
\end{align}
Since
\begin{align*}
	I
		&= \grad \x
		= \grad_\x(\eta(t_1, t_2; \eta(t_2, t_1; \x)))
		= \grad \eta(t_1, t_2; \eta(t_2, t_1; \x)) \grad \eta(t_2, t_1; \x),
\end{align*}
we have
\begin{align*}
	\grad \eta(t_1, t_2; \eta(t_2, t_1; \x))
		= (\grad \eta(t_2, t_1; \x))^{-1}.
\end{align*}
Hence,
\begin{align}\label{e:BInv}
	B_+^{-1}
		= \grad \eta(t, \time; \x), \quad
	B_-^{-1}
		= \grad \eta(t, 0; \x).
\end{align}

In what follows, we will apply the chain rule quite liberally, and will avoid much redundancy by largely employing the time-space derivative operator $D$ of \cref{e:D}.

Since $\time = 0$ and $\pos = \pos_0$ on $S$, we can write \cref{e:Dpos} as
\begin{align}\label{e:Dpospos0}
	D \pos
		&= D \pos_0 + \uu_0(\pos_0) \otimes D \time
			\onS.
\end{align}
On $\olUPlus$, we can see from the chain rule that
\begin{align*}
	\begin{split}
	\prt_t M
		&= \prt_t (\grad \eta(t, 0; \x)) - \prt_t (\grad \eta(t, \time; \x)) \\
		&= \prt_{t_1} \grad \eta(t, 0; \x) - \prt_{t_1} \grad \eta(t, \time; \x)
			-\prt_{t_2} \grad \eta(t, \time; \x) \prt_t \time, \\
	\grad M
		&= \grad_\x (\grad \eta(t, 0; \x)) - \grad_\x (\grad \eta(t, \time; \x)) \\
		&= \grad \grad \eta(t, 0; \x) - \prt_{t_2} \grad \eta(t, \time; \x) \otimes \grad \time
			- \grad \grad \eta(t, \time; \x).
	\end{split}
\end{align*}
But, from \cref{L:prt2gradeta}, on $S$, where $\time = 0$,
$
	\prt_{t_2} \grad \eta(t, \time; \x)
		= \grad \uu(0, \eta(t, 0; \x)) \grad \eta(t, 0; \x)
		= \grad \uu_0(\pos_0) B_-^{-1},
$
so we have,
\begin{align}\label{e:prtgradM}
	\begin{split}
	\prt_t M
		&= -\prt_t \time \grad \uu_0(\pos_0) B_-^{-1}, \quad
	\grad M
		= - \grad \uu_0(\pos_0) B_-^{-1} \otimes \grad \time, \\
		D M
			&= - \grad \uu_0(\pos_0) B_-^{-1} \otimes D \time
			\onS.
	\end{split}
\end{align}

Letting $\ZZ_\pm := \Y_\pm - \G_\pm$, we have, from \cref{e:LagrangianForm},
\begin{align}\label{e:YminusYplusAlt}
	\begin{split}
		B_-^{-1} \ZZ_-(t, \x)
			= \Y_0(\pos_0), \quad
		B_+^{-1} \ZZ_+(t, \x)
			&= \H(\time, \pos).
	\end{split}
\end{align}
Using the expressions for $\prt_t G$ and $\grad G$ in \cref{P:Forcing},
\begin{align*}
	\begin{split}
	\prt_t (\G_- - \G_+)(t, \x)
		&= (\eta(0, t)_* \g(0))(t, \x) \prt_t \time(t, \x)
		= B_- \g(0, \pos_0) \prt_t \time(t, \x)
			\onS, \\
	\grad (\G_- - \G_+)(t, \x)
		&= (\eta(0, t)_* \g(0))(t, \x) \otimes \grad \time
		= B_- \g(0, \pos_0) \otimes \grad \time
			\onS,
	\end{split}
\end{align*}
or,
\begin{align}\label{e:DGG}
	D (\G_- - \G_+)(t, \x)
		&= B_- \g(0, \pos_0) \otimes D \time
			\onS.
\end{align}
Observe that by virtue of \cref{P:Forcing}, our assumptions on $\g$ were enough to obtain \cref{e:DGG}.

\begin{proof}[\textbf{Proof of \cref{P:N1Condition}}]
	Because $B_-^{-1}$ is invertible,
	$D \Y_- = D \Y_+ \onS \text{ iff } B_-^{-1} D (\Y_+ - \Y_-) = 0$. Then,
	using \cref{e:YminusYplusAlt,e:DGG},
	\begin{align*}
		B_-^{-1} D &(\Y_+ - \Y_-)
			= B_-^{-1} D (\ZZ_+ - \ZZ_-)
				+ B_-^{-1} D (\G_+ - \G_-) \\
			&= B_-^{-1} D (\ZZ_+ - \ZZ_-)
				- (B_-^{-1} B_-) \g(0, \pos_0) \otimes D \time
			= B_-^{-1} D (\ZZ_+ - \ZZ_-)
				- \g(0, \pos_0) \otimes D \time
	\end{align*}
	\onS. Also on $S$,
	\begin{align*}
		B_-^{-1} D (\ZZ_+ - \ZZ_-)
			&=  D (B_-^{-1} (\ZZ_+ - \ZZ_-))
				- (D B_-^{-1}) (\ZZ_+ - \ZZ_-)
			= D (B_-^{-1} (\ZZ_+ - \ZZ_-)),
	\end{align*}
	since $\ZZ_+ - \ZZ_- = \H(\time, \pos) - \Y_0(\pos_0) = 0$ on $S$ by $\cond_0$.
	So,
	\begin{align}\label{e:BminusinvZZ}
		\begin{split}
		B_-^{-1} D (\ZZ_+ - \ZZ_-)
			&= D (B_+^{-1} \ZZ_+ - B_-^{-1} \ZZ_-)
				+ D ((B_-^{-1} - B_+^{-1}) \ZZ_+)
			\\
			&= D [\H(\time, \pos) - \Y_0(\pos_0)]
				+ D (M \ZZ_+)
				\onS.
		\end{split}
	\end{align}
	
	Then,
	\begin{align*}
		D [\H(&\time, \pos) - \Y_0(\pos_0)]
			= \prt_{t_1} \H(\time, \pos) \otimes D \time
				+ \grad_\Gamma \H(\time, \pos) D \pos
				-  \grad \Y_0(\pos_0) D \pos_0
	\end{align*}
	on $S$. But, using \cref{e:Dpospos0},
	\begin{align*}
		\grad_\Gamma \H(\time, &\pos) D \pos
				-  \grad \Y_0(\pos_0) D \pos_0
			= [\grad_\Gamma \H(\time, \pos) - \grad \Y_0(\pos_0)] D \pos
				+ \grad \Y_0(\pos_0) [D \pos - D \pos_0] \\
			&=  [\grad_\Gamma \H(\time, \pos) - \grad \Y_0(\pos_0)] D \pos
				+ \grad \Y_0(\pos_0)
					(\uu_0(\pos_0) \otimes D \time)
				\onS.
	\end{align*}

	Now, $\pos = \pos(t, \x)$ always lies on $\Gamma_+$, so
	\begin{align*}
		[\grad_\Gamma \H(\time, \pos) - \grad \Y_0(\pos_0)] D \pos
			= 0,
	\end{align*}
	where we used $\cond_0$ along with $N = 1$ regularity.
	Hence, \onS,
	\begin{align}\label{e:DHminusY0}
		D [\H(&\time, \pos) - \Y_0(\pos_0)]
			= [\prt_{t_1} \H(\time, \pos)
				+ \grad \Y_0(\pos_0) \uu_0(\pos_0)]
				\otimes D \time.
	\end{align}
	
	And, using \cref{e:prtgradM}, the vanishing of $M$ on $S$, and $\cond_0$ once more,
	\begin{align*}
		D (M \ZZ_+)
			&= M D \ZZ_+ + D M \ZZ_+
			= - [\grad \uu_0(\pos_0) B_-^{-1} \otimes D \time] \ZZ_+ \\
			&= - \grad \uu_0(\pos_0) B_-^{-1} \ZZ_- \otimes D \time
			= - \grad \uu_0(\pos_0) \Y_0(\pos_0) \otimes D \time
			\onS.
	\end{align*}
	
	Combined, these calculations give
	\begin{align*}
		B_-^{-1} D (\Y_+ - \Y_-)
			&= [\prt_{t_1} \H(\time, \pos)
				+ \grad \Y_0(\pos_0) \uu_0(\pos_0) - \grad \uu_0(\pos_0) \Y_0(\pos_0)
					- \g(0, \pos_0)]
				\otimes D \time
	\end{align*}
	on $S$. Recalling that $B_-$ is always invertible,
	if $\cond_1$ holds then the right-hand side vanishes so $D (\Y_+ - \Y_-) = 0 \onS$. Conversely,
	if the left-hand side vanishes then since $\prt_t \time > 0$ up to at least time $T^*$,
	we know that $D \time$ does not vanish up to $T^*$, so $\cond_1$ must hold. 
\end{proof}

\section{Vorticity}\label{S:Vorticity}

\noindent
\begin{prop}\label{P:YInRangeOfCurl}
	Assume that for some $N \ge 0$ the data has regularity $N$ as in \cref{D:NReg},
	$\cond_N$ holds, $\Y_0$ and $\g(t)$ for all $t \in [0, T]$ are in the range of the curl,
	and \cref{e:HigherH,e:RangeOfCurlCond} hold.
	Then the solution $\Y$ to \Our given by
	\cref{P:LinearExistenceFirstPart}
	(unique classical for $N \ge 1$, Lagrangian for all $N \ge 0$)
	is in the range of the curl for the lifetime of the solution.
\end{prop}

\begin{remark}\label{R:NoPenetrationBCs}
	\cref{P:YInRangeOfCurl} along with \cref{P:LinearExistenceFirstPart}
	(and \cref{R:BeyondTStar})
	give a complete proof of \cref{T:LinearExistence}.
	We note that
	if $\Gamma_0 = \Gamma$, the classical setting of impermeable boundary
	conditions, no vorticity is transported off of the boundary,
	so $U_-$ is all of $Q$ and many
	of the flow map constructs, such as $S$, $\time$, and $\pos$, are unnecessary.
\end{remark}

In this section, we make use of the $\dv_\Gamma$ operator,
so let us introduce the few facts we will need about it, referring the interested reader to \NonlinearPaper{Appendix B} for more details (including its explicit expression in coordinates on the boundary). First, we define $\grad_\Gamma$, the gradient in the tangent space to $\Gamma$, by requiring that for any $f \in C^1(\Gamma)$ and any continuously differentiable curve $\x(s)$ on $\Gamma$ parameterized by arc length,
\begin{align*}
	\grad_\Gamma f \cdot \x'(0)
		&= \lim_{s \to 0} \frac{f(\x(s)) - f(\x(0))}{s}.
\end{align*}
We then define $\dv_\Gamma$ by duality, requiring that for any $f \in C^1(\Gamma)$, $\vv \in C^1(\Gamma)^d$,
\begin{align}\label{e:IBP}
	\int_\Gamma \vv \cdot \grad_\Gamma f
		&= - \int_\Gamma \dv_\Gamma \vv \, f.
\end{align}
See, for instance, Proposition 2.2 of \cite{T1996I}.

We also need the following two identities (again, see \NonlinearPaper{Appendix B}):
For all $f \in C^1(\Gamma)$,
\begin{align}\label{e:dvfv}
	\dv_\Gamma (f \vv)
		&= f \dv_\Gamma \vv + \vv \cdot \grad_\Gamma f
\end{align}
and for any $C^1$ vector field $\vv$ defined in an epsilon neighborhood $\Cal{N}$ of $\Gamma$,
\begin{align}\label{e:divdivGamma}
	\dv \vv = \dv_\Gamma \vv|_\Gamma + \grad (\vv \cdot \n) \cdot \n + \kappa \vv \cdot \n,
\end{align}
where $\kappa = \kappa_1 + \kappa_2$ is the mean curvature and $\n$ is extended orthogonally into $\Cal{N}$.

\begin{proof}[\textbf{Proof of \cref{P:YInRangeOfCurl}}]
From \cref{T:VectorPotential}, we know that $\Y(t)$ will be in the range of the curl if and only if it is divergence-free and has vanishing external flux through each boundary component, as defined in \cref{e:ExternalFlux}.

That $\dv \Y_-$ is transported on $U_-$ is essentially classical. It comes from taking the divergence of \OurOne, which yields, after a few calculations, that (since $\dv \g = \dv \curl \f = 0$),
\begin{align}\label{e:divTransported}
	\prt_t \dv \Y + \uu \cdot \grad \dv \Y = 0.
\end{align}
This holds in weak form when $N = 0$, but still gives $\dv \Y = 0$ on $\ol{U}_-$.

To show that $\dv \Y = 0$ on $\ol{U}_+$, we need only show that $\dv \Y = 0$ on $\Gamma_+$, because \cref{e:divTransported} holds on all of $\ol{Q}$
and all flow lines on $\ol{U}_+$ originate on $[0, T] \times \Gamma_+$. Complicating things, however, is that values of $\Y$ are generated on $\Gamma_+$ by $\H$.

Let us start by assuming that $N \ge 1$. Then because
\begin{align*}
	\dv [\Y \otimes \uu - \uu \otimes \Y]^i
		&= \prt_j [Y^i u^j - u^i Y^j]
		= Y^i \dv \uu + [\uu \cdot \grad \Y]^i
			- u^i \dv \Y - [\Y \cdot \grad \uu]^i,
\end{align*}
\OurOne gives
\begin{align}\label{e:YdivCalcOnQ}
	\prt_t \Y
		+ \dv [\Y \otimes \uu - \uu \otimes \Y]
		+ (\dv \Y) \, \uu
		= \g.
\end{align}
Restricting \cref{e:YdivCalcOnQ} to $[0, T] \times \Gamma_+$ and taking the inner product with $\n$, we have,
\begin{align}\label{e:PreHCond}
	U^\n \dv \Y
		&= \g \cdot \n - \prt_t \Y \cdot \n  - \dv [\Y \otimes \uu - \uu \otimes \Y] \cdot \n 
			\text{ on } [0, T] \times \Gamma_+.
\end{align}
Using \cref{L:divdiv}, below, the condition for $\dv \Y = 0$ on $[0, T] \times \Gamma_+$ and hence on all of $Q$ can be written (noting that $U^\n$ never vanishes on $[0, T] \times \Gamma_+$),
\begin{align}\label{e:HCond}
	\prt_t H^\n
		+ \dv_\Gamma [H^\n \uu^\BoldTau  - U^\n \H^\BoldTau]
		- \g \cdot \n
		= 0
			\text{ on } [0, T] \times \Gamma_+,
\end{align}
which is \cref{e:RangeOfCurlCond}.
We conclude that $\dv \Y = 0$ on $\ol{U}_+$ and hence on $Q$. (Observe that $\dv \Y_0 = 0$ on $\ol{\Omega}$ was important to conclude that $\dv \Y(t) = 0$ on $S(t)$.)

It follows from \cref{e:HCond} that the external flux through any component $\Gamma_\ell$ of $\Gamma_+$ vanishes, since $\dv_\Gamma [H^\n \uu^\BoldTau  - U^\n \H^\BoldTau]$ integrates to zero on $\Gamma_\ell$, as does $\g \cdot \n$ by \cref{T:VectorPotential}. Then, because $\prt_t Y^n = \prt_t H^\n \in C^\al([0, T] \times \Gamma_+)$, even for $N = 0$, we have
\begin{align*}
	\diff{}{t} \Phi_\ell^\Gamma(\Y(t))
		&= \diff{}{t} \int_{\Gamma_\ell}
			\brac{\g \cdot \n - \dv_\Gamma [H^\n \uu^\BoldTau  - U^\n \H^\BoldTau]}
		= 0.
\end{align*}
Hence, $\Phi_\ell^\Gamma(\Y) = 0$, since $\Y_0$ is in the range of the curl.

For $N \ge 1$, we have sufficient regularity to obtain the vanishing of the external fluxes through any boundary component $\Gamma_\ell$. Returning to \cref{e:YdivCalcOnQ}, using that $\dv \Y = 0$ and \cref{L:divdiv}, we have, for any boundary component $\Gamma_\ell$,
\begin{align*}
	\prt_t Y^n + \dv_\Gamma [Y^n \otimes \uu^\BoldTau - U^\n \otimes \Y^\BoldTau]
			- \g \cdot \n = 0
		\text{ on } [0, T] \times \Gamma_\ell.
\end{align*}
Now we have sufficient regularity to make the same calculation we made on components of $\Gamma_+$ to obtain $\Phi_\ell^\Gamma(\Y) = 0$.

Now suppose that $N = 0$. We will make an approximation argument like that in the proof of \cref{L:RankineHuginoit}, but without approximating $\Y_0$.

Let $(\uu_i)$ be the sequence in $\uBoundarySpaceOne{U_i^\n}$ with $\uu_i \to \uu$ in $\uInputSpace$ given by \cref{L:VelDense01}. Let $\Y_i$ be the Lagrangian solution to \Our as in \cref{D:LagrangianSolution} with $\uu_i$ in place of $\uu$, and let $U_{i, \pm}$ be the $U_\pm$ sets corresponding to $\uu_i$. Because we did not change $\Y_0$, $\cond_0$ is satisfied, though $\cond_1$ need not be, so we only have that $\Y_i \in C^{1, \al}(U_{i, +})$ and, from \cref{P:RegNFromS} applied for $N = 0$, that $\Y_i \in C^\al(Q)$. We have sufficient regularity of $\Y_i$ on $U_{i, +}$, however, to conclude, after applying \cref{L:divdiv} and using that $\Y_i = \H$ on $[0, T] \times \Gamma_+$, that
\begin{align*}
	\prt_t H^\n
		+ \dv_\Gamma [H^\n \uu_i^\BoldTau  - U_i^\n \H^\BoldTau]
		+ (\dv \Y_i) U_i^\n
		= \g \cdot \n,
\end{align*}
holds classically for all $t > 0$. It follows that
\begin{align}\label{e:dvYnUn}
	\begin{split}
	(\dv \Y_i) U_i^\n
		&= \g \cdot \n -\prt_t H^\n
			- \dv_\Gamma [H^\n \uu_i^\BoldTau  - U_i^\n \H^\BoldTau] \\
		&= \g \cdot \n -\prt_t H^\n
			- \dv_\Gamma [H^\n \uu^\BoldTau  - U^\n \H^\BoldTau]
			+ \dv_\Gamma [(U_i^\n - U^\n) \H^\BoldTau] \\
			&\qquad
			+ \dv_\Gamma [H^\n (\uu^\BoldTau - \uu_i^\BoldTau)] \\
		&= \dv_\Gamma [(U_i^\n - U^\n) \H^\BoldTau] 
			+ \dv_\Gamma [H^\n (\uu^\BoldTau - \uu_i^\BoldTau)] \\
		&= (U_i^\n - U^\n) \dv_\Gamma \H^\BoldTau
			+ \grad_\Gamma (U_i^\n - U^\n) \cdot \H^\BoldTau
			\\
			&\qquad
			+ H^\n \dv_\Gamma (\uu^\BoldTau - \uu_i^\BoldTau)
			+ \grad_\Gamma H^\n \cdot (\uu^\BoldTau - \uu_i^\BoldTau),
	\end{split}
\end{align}
where we used the assumption that \cref{e:RangeOfCurlCond} holds. But $\uu_i \to \uu$ in $\uInputSpace$ so, in particular, $\uu_i \to \uu$ in $C^\al([0, T]; C^{1, \al}(\Gamma_+))$. Hence, taking advantage of \cref{e:HigherH}, $\dv \Y_i \to 0$ in $C^\al(\Gamma_+)$, and arguing as above using transport, $\norm{\dv \Y_i}_{C^\al(U_{i, +})} \to 0$.

Let $\varphi \in \Cal{D}(U_+)$. The convergence $\uu_i \to \uu$ in $\uInputSpace$ gives that $\supp \varphi \subseteq U_{i, +}$ for all sufficiently large $i$. We can see, then, from the form of $\Y_\pm$ and $\Y_{i, \pm}$ coming from \cref{e:LagrangianForm}, that $\Y_i \to \Y$ in $C^\al(\supp \varphi)$. Then,
\begin{align*}
	(\dv (\Y_i - \Y), \varphi)
		&= - (\Y_i - \Y, \grad \varphi)
		\to 0.
\end{align*}
But also $(\dv \Y_i, \varphi) \to 0$ since $\norm{\dv \Y_i}_{C^\al(U_{i, +})} \to 0$, so it follows that as a distribution, $\dv \Y = 0$ on $U_+$.

The vanishing of the external fluxes on components of $\Gamma_+$ holds directly for $N = 0$ as noted above. Finally, the form of $\Y_\pm$ and $\Y_{i, \pm}$ coming from \cref{e:LagrangianForm} also shows that $\Y_i \to \Y$ in $C^\al(V)$ for any $V$ compactly supported in $U_- \cup ([0, T] \times (\Gamma_- \cup \Gamma_0))$, from which the vanishing of
all external fluxes of $\Y$ follows, so, in fact, $\Y(t)$ remains in the range of the curl.
\end{proof}

\begin{remark}\label{R:HHigherRegularity}
	In the proof above we used that $\H \in C^{1, \al}([0, T] \times \Gamma_+)$ and
	$\g \in C^\al(Q)$ for both $N = 0$ and $N = 1$, so that we could make the approximation
	argument by only approximating $\uu$ with a higher-regularity sequence.
	If for $N = 0$ we assumed, as would seem natural, only that
	$\H \in C^\al([0, T] \times \Gamma_+)$, we would also need to approximate $\H$
	by some $\H_i$. But because of the factor $\grad_\Gamma H^\n$ in \cref{e:dvYnUn},
	it would be  difficult, if not impossible, to successfully complete the
	approximation argument to conclude that $\dv \Y = 0$ when $N = 0$.
\end{remark}

We used \cref{L:divdiv} in the proof of \cref{P:YInRangeOfCurl}.
\begin{lemma}\label{L:divdiv}
	Assume that $N \ge 1$. We have
	\begin{align}\label{e:dvYuuYn}
		\dv [\Y \otimes \uu - \uu \otimes \Y] \cdot \n
			&= \dv_\Gamma [Y^\n \uu^\BoldTau  - U^\n \Y^\BoldTau]
	\end{align}
	with
	\begin{align*}
		Y^\n \uu^\BoldTau - U^\n \Y^\BoldTau
			= H^\n \uu^\BoldTau - U^\n \H^\BoldTau
			\text{ on }\Gamma_+.
	\end{align*}
	On any boundary component $\Gamma_\ell$,
	\begin{align}\label{e:IntdvYuuYnZero}
		\int_{\Gamma_\ell}
				\dv [\Y \otimes \uu - \uu \otimes \Y] \cdot \n
			=
			\int_{\Gamma_\ell}
				\dv_\Gamma [Y^\n \uu^\BoldTau  - U^\n \Y^\BoldTau]
				= 0.
	\end{align}
\end{lemma}
\begin{proof}
	Let $M := \Y \otimes \uu - \uu \otimes \Y$, an anti-symmetric matrix-valued function.
	Then, working in Cartesian coordinates, 
	\begin{align*}
		\dv M \cdot \n
			&= \prt_k M_k^i n^i
			= \prt_k (M_k^i n^i) - M_k^i \prt_k n^i
			= \dv (M \cdot \n) - M_k^i \prt_k n^i
			= \dv (M \cdot \n).
	\end{align*}
	Here, $M_k^i \prt_k n^i = 0$ because $M$ is antisymmetric while $\grad \n$ is symmetric
	(since, as in the proof of \cref{L:Surface}, we can write locally,
	$\Gamma = \varphi^{-1}(0)$ for some $\varphi \in C^{1, \al}(\R^3)$,
	which we can choose so that $\grad \varphi|_{\varphi_0^{-1}(0)} = \n$
	on $\Gamma$).
	By \cref{e:divdivGamma},
	\begin{align*}
		\dv (M \cdot \n)
			= \dv_\Gamma (M \cdot \n)|_\Gamma
				+ \grad ((M \cdot \n) \cdot \n) \cdot \n
				+ \kappa (M \cdot \n) \cdot \n
			= \dv_\Gamma (M \cdot \n)|_\Gamma,
	\end{align*}
	since $(M \cdot \n) \cdot \n = M_i^j n^i n^j = - M_j^i n^j n^i = - M_i^j n^i n^j$,
	so $\grad ((M \cdot \n) \cdot \n) \cdot \n = \kappa (M \cdot \n) \cdot \n = 0$.
	
	Finally, $M \cdot \n = Y^\n \uu - U^\n \Y$, and we obtain \cref{e:dvYuuYn}.
	
	Integrating by parts along the boundary using \cref{e:IBP} gives
	\cref{e:IntdvYuuYnZero}, completing the proof.
\end{proof}

\section{Velocity}\label{S:Velocity}

\noindent In this section, we prove \cref{T:VelocityForm}. We will show that the operator $\VV_c$ of \cref{L:Vc} recovers the harmonic component of the velocity field corresponding to a solution to \Our.

For a vorticity field $\bmu \in \vortSpace$, define the matrix-valued function,
\begin{align*}
	\bOmega(\bmu) :=  \grad K[\bmu] - (\grad K[\bmu])^T,
\end{align*}
noting that the nonzero components of $\bOmega(\bmu)$ are those of $\pm \bmu$.

\begin{lemma}\label{L:Vc}
	Assume the data has regularity $N \ge 0$ and let
	$K$ be the Biot-Savart operator of \cref{T:VectorPotential}.
	For a velocity field $\uu$ and vorticity field $\bmu$, define the vector field,
	\begin{align*}
		(\VV_c (\uu, \bmu))(t)
			:= P_{H_c} &\uu_0 + \int_0^t P_{H_c} \f(s) \, ds
				- \int_0^t P_{H_c} P_H \pr{\uu(s) \cdot \bOmega(\bmu)(s)} \, ds.
	\end{align*}
	Then $\VV_c \colon \uInputSpaceN \times \vortSpaceN \to C^{N + 1, \al}(Q) \cap C([0, T]; H_c)$.
\end{lemma}
\begin{proof}
	From the regularity of $\uu$ and $\bmu$, $\uu(s) \cdot \bOmega(\bmu)(s) \in C^{N, \al}(Q)$.
	By \cref{L:LerayContinuity,L:TamingHc}, it follows that
	$P_{H_c} P_H \pr{\uu(s) \cdot \bOmega(\bmu)(s)} \in \uInputSpaceN$.
	The integrals in time along with the assumed regularity of $\uu_0$
	and $\f$ then give that $\VV_c (\uu, \bmu) \in \uSolSpaceN$.
\end{proof}

\begin{proof}[\textbf{Proof of \cref{T:VelocityForm}}]
	Letting $\bomega = \Y$, we write \Our as
	\begin{align}\label{e:OurOmega}
		\begin{cases}
			\prt_t \bomega + \dv (\bomega \otimes \uu)  - \bomega \cdot \grad \uu = \g
					&\text{in } Q, \\
			\bomega = \H
				&\text{on } [0, T] \times \Gamma_+, \\
			\bomega(0) = \bomega_0
				&\text{on } \Omega.
		\end{cases}
	\end{align}
	Now, \cref{e:OurOmega}$_1$ holds at least in the sense of elements of $\Cal{D}'(Q)$.
	We know from \cref{T:LinearExistence} that $\bomega$ is in the range of the curl,
	so from \cref{e:KUDef}, we know that $\bomega = \curl \vv$, where
	\begin{align*}
		\vv
			= K_{U^\n}[\bomega] + \vv_c
			= \ol{\vv} + \vv_c + \VV, \quad
		\ol{\vv}(t) := K[\bomega(t)] \in H_0,
	\end{align*}
	where we defer the free choice of $\vv_c(t) \in H_c$ until later.

	A direct calculation gives
	\begin{align*}
		\dv (\bomega \otimes \uu) - \bomega \cdot \grad \uu
			= \curl(\uu \cdot \grad \vv - \uu \cdot (\grad \vv)^T).
	\end{align*}
	Thus,
	\begin{align*}
		\curl \pr{\prt_t \vv + \uu \cdot \grad \vv - \uu \cdot (\grad \vv)^T - \f}
			&= 0,
	\end{align*}
	equality holding in $\Cal{D}'(Q)$.
	We conclude, since $\bomega = \curl \ol{\vv} = \curl \vv$, that
	\begin{align}\label{e:NSwithw}
		\prt_t \vv + \uu \cdot \bOmega(\bomega) - \f
			&= \prt_t \vv + \uu \cdot \grad \vv - \uu \cdot (\grad \vv)^T - \f
			= -\grad \pi + \zz,
	\end{align}
	where $\zz(t) \in H_c$ and $\grad \pi(t) \in L^2(\Omega)$.
	Since $\VV$ is a gradient,
	$P_H \VV = 0$, so we see that $P_{H_0} P_H \vv = \ol{\vv}$.
	Letting $\ol{\f} = P_{H_0} \f$, $\f_c = P_{H_c}\f$, we first apply $P_H$ to
	both sides of \cref{e:NSwithw}, giving
	\begin{align}\label{e:NSH0}
		\prt_t (\ol{\vv} + \vv_c)
			+ P_H \pr{\uu \cdot \bOmega(\bomega) }
			- \f
			&= \zz
	\end{align}
	then apply $P_{H_0}$, $P_{H_c}$ to give
	\begin{align}\label{e:prtvvc}
		\begin{split}
		\prt_t \ol{\vv}
			+ P_{H_0} P_H \pr{\uu \cdot \bOmega(\bomega)}
			- \ol{\f} &= 0, \\
		\prt_t \vv_c
			+ P_{H_c} P_H \pr{\uu \cdot \bOmega(\bomega)}
			- \f_c
			&= \zz.
		\end{split}
	\end{align}
	We now choose $\vv_c(t) = \VV_c (\uu, \bomega)(t)$,
	from which $\zz = 0$ follows,
	and	\cref{e:NSwithw} becomes \cref{e:LinearVelEq}.
	
	Now suppose that $\bomega = \curl \uu$. Then also $\curl \uu = \curl \vv$,
	so $\vv = \uu + \ww$ for some $\ww(t) \in H_c \cap \uInputSpaceHomoN$,
	noting that $\uInputSpaceHomoN$ of \cref{e:BoundarySpace}
	is defined as $\uInputSpaceN$, but with
	$\uu \cdot \n = 0$ on $\prt \Omega$, and
	where we appealed to \cref{L:KRegular}. Also,
	$
		\uu \cdot \bOmega(\bomega) 
			= \uu \cdot \grad \uu - \uu \cdot (\grad \uu)^T
	$,
	since the off-diagonal components of the antisymmetric $\grad \ww - (\grad \ww)^T$ come
	from the components of $\curl \ww = 0$. Hence, \cref{e:LinearVelEq} becomes
	\begin{align*}
		\prt_t \uu + \uu \cdot \grad \uu - \uu \cdot (\grad \uu)^T + \grad \pi
			= \prt_t \uu + \uu \cdot \grad \uu + \grad p
			= \f + \prt_t \ww,
	\end{align*}
	where $p = \pi - (1/2) \abs{\uu}^2$. If, further, $\vv = \uu$ then $\prt_t \ww = 0$,
	and we recover \cref{e:EulerInflowOutflow} with $\uu \cdot \n = U^\n$
	on $[0, T] \times \Gamma$.
	
	From \cref{L:Vc}, $\vv_c \in C^{N + 1, \al}(Q) \cap C([0, T]; H_c)$. From \cref{L:KRegular},
	$\ol{\vv} = K_{U^\n}[\bomega] \in \uInputSpaceN$. And from \cref{e:prtvvc}$_1$ and
	\cref{L:LerayContinuity}, $\prt_t \ol{\vv} \in C^{N, \al}(Q)$, which gives the
	additional time continuity to conclude that $\vv \in \uSolSpaceN$.
	Returning to \cref{e:NSwithw} (where now $\zz = 0$), we conclude that
	$\grad \pi \in C^{N, \al}(Q)$.

To prove uniqueness, suppose that $\prt_t \vv_j + \uu \cdot \grad \vv_j - \uu \cdot (\grad \vv_j)^T + \grad \pi_j = \f$ for $j = 1, 2$. Letting $\vv = \vv_1 - \vv_2$, we have $\prt_t \vv + \uu \cdot (\grad \vv - (\grad \vv)^T) + \grad \pi = 0$, where $\pi = \pi_1 - \pi_2$. But also, $\curl \vv = \bomega - \bomega = 0$ so $\grad \vv - (\grad \vv)^T = 0$, and we see that
$
	\prt_t \vv = -\grad \pi
$.
Then $\vv \cdot \n = 0$ so $\vv \in H$, meaning that it must be that $\grad \pi = 0$ and hence $\prt_t \vv = 0$. Finally, since $\vv(0) = \uu_0 - \uu_0 = 0$, we conclude that $\vv = 0$ on $Q$, giving uniqueness.
\end{proof}

\section*{Acknowledgements}

\noindent
The authors thank an anonymous referee for valuable comments and suggestions that improved the exposition of this paper.
Gie was partially supported by a Simons Foundation Collaboration Grant for Mathematicians; Research R-II Grant and the Ascending Star Fellowship, Office of EVPRI, University of Louisville; Brain Pool Program through the National Research Foundation of Korea (NRF) (grant number: 2020H1D3A2A01110658).
Mazzucato was partially supported by the US National Science Foundation Grant DMS-1909103. Part of this work was prepared while Kelliher and Mazzucato were participating in a program hosted by the Mathematical Sciences Research Institute in Berkeley, California, in Spring 2021, supported by the National Science Foundation under Grant No. DMS-1928930.
Mazzucato would like to thank the Isaac Newton Institute for
Mathematical Sciences, Cambridge, for support and hospitality during the
programme, \textit{Mathematical aspects of turbulence: where do we stand?},
where work on this paper was partially undertaken.
The work of the Institute is supported by EPSRC grant no EP/R014604/1.

\appendix
\section{\Holder space lemmas}\label{S:HolderLemmas}

\noindent
We collect here a number of estimates in \Holder spaces.

For any $V \subseteq \R^d$, $d \ge 1$, define the classical \Holder space, $C^\al(V)$, with the norm
\begin{align*}
	\norm{f}_{C^\al(V)}
		:= \norm{f}_{L^\iny(V)}
			+ \sup_{\y_1 \ne \y_2 \in V}
			\frac{\abs{f(\y_1) - f(\y_2)}}
				{\abs{\y_1 - \y_2}^\al}.
\end{align*}

\begin{lemma}\label{L:HolderProd}
	Let $f, g \in C^\al(U)$. Then
	\begin{align*}
		\norm{fg}_{C^\al}
			&\le \norm{f}_{C^\al} \norm{g}_{C^\al}, \\
		\norm{fg}_{\dot{C}^\al}
			&\le \norm{f}_{L^\iny} \norm{g}_{\dot{C}^\al}
				+ \norm{g}_{L^\iny} \norm{f}_{\dot{C}^\al}, \\
		\norm{fg}_{C^\al}
			&\le \norm{f}_{L^\iny} \norm{g}_{L^\iny} 
				+ \norm{f}_{L^\iny} \norm{g}_{\dot{C}^\al}
				+ \norm{g}_{L^\iny} \norm{f}_{\dot{C}^\al} \\
			&\le
				\norm{f}_{L^\iny} \norm{g}_{C^\al}
				+ \norm{g}_{L^\iny} \norm{f}_{C^\al}.
	\end{align*}
\end{lemma}
\begin{proof}
	These are all classical.
\end{proof}

\begin{lemma}\label{L:HolderComp}
	Let $U, V$ be open subsets of Euclidean spaces, $\al \in (0, 1]$, and $k \ge 1$ an integer. 
	If $f \in C^{k, \al}(U)$ and $g \in C^{k + 1, \al}(V)$ with $g(V) \subseteq U$
	then
	\begin{align}\label{e:HolderComp}
		\begin{split}
		\norm{f \circ g}_{\dot{C}^\al(V)}
			&\le \norm{f}_{\dot{C}^\al(U)} \norm{g}_{Lip(V)}^\al, \\
		\norm{f \circ g}_{C^\al(V)}
			&\le \norm{f}_{L^\iny(U)}
				+ \norm{f}_{\dot{C}^\al(U)} \norm{g}_{Lip(V)}^\al
			\le \norm{f}_{C^\al(U)} \brac{1 + \norm{g}_{Lip(V)}^\al}, \\
		\norm{f \circ g}_{C^{k, \al}(V)}
			&\le C(k) \norm{f}_{C^{k, \al}(U)} \brac{1 + \norm{g}_{C^{k + 1}(V)}}^{k + 1},
		\end{split}
	\end{align}
	where $Lip$ is the homogeneous Lipschitz semi-norm
	and $\dot{C}^\al$ is the homogeneous \Holder norm.
\end{lemma}
\begin{proof}
	These bounds are all classical.
\end{proof}

\section{The continuity of the Biot-Savart law}\label{A:KTechnical}

\noindent The operator $K$ of \cref{T:VectorPotential} allows us to recover a divergence-free vector field in $H_0 \cap H^1(\Omega)^3$ from its curl, $\bomega$. We need, however, to obtain estimates on $K[\bomega]$ in terms of $\bomega$ in various norms. To do that, we will use results from Kato, Mitrea, Ponce, and Taylor's \cite{KMPT2000}. For this, we need to explore the Hodge decomposition slightly further than we did in \cref{S:BSLaw}.

Let $\Cal{O}_j$ be the component of $\R^3 \setminus \ol{\Omega}$ whose boundary is $\Gamma_j$, $j = 1, \dots, b + 1$. 
Let $\Sigma_1, \dots, \Sigma_M$ be pairwise disjoint $C^{N, \al}$-regular surfaces (``admissible cuts'') which, when removed from $\Omega$ render it simply connected.
Let $\vv$ lie in the space $H$ of \cref{e:HSpace}, so $\vv \cdot \n \in H^{-\frac{1}{2}}(\Sigma)$. The internal flux $\Phi_i$ of $\vv$ across $\Sigma_i$ is defined to be the value of
\begin{align}\label{e:Flux}
	\Phi_i(\vv)
		:= \int_{\Sigma_i} \vv \cdot \n,
\end{align}
where the direction of the unit normal vector $\n$ to $\Sigma_i$ is fixed by an arbitrarily chosen orientation to $\Sigma_i$. Because $\vv$ is divergence-free and tangential to the boundary, it is easy to see that the internal fluxes do not depend upon the specific choices of the $\Sigma_i$. It is classical (going back in some form to Helmholtz) that
\begin{align*}
	H_0
		= \set{\vv \in H \colon \text{all internal fluxes are zero}}.
\end{align*}

Fix, arbitrarily, points $y_j \in \Cal{O}_j$ for each $j = 1, \dots, b + 1$ and define
\begin{align*}
	g_j(x)
		:= \grad G(\cdot - y_j),
\end{align*}
where $G(x) := -1/(4 \pi\abs{x})$ is the fundamental solution of the Laplacian. Note, then, that $\dv g_j$ and $\curl g_j$ both vanish away from $y_j$.

\begin{theorem}\label{T:KMPT}
	Assume that $\Gamma$ is $C^{k + 1, \al}$-regular, $k \ge 0$, and
	let $\bomega \in C^{k, \al}(\Omega)$ (or $\bomega \in H^{k, p}(\Omega)$, $p \in (1, \iny)$).
	There exists an antisymmetric matrix-valued
	function $M \in C^{k + 1, \al}(\Omega)$ (or $M \in H^{k + 1, p}(\Omega)$) such that
	\begin{align*}
		\bomega = \dv M + \sum_j \la_j g_j,
	\end{align*}
	where $\dv M$ is the row-by-row divergence of the matrix $M$
	(observe that $\dv \dv M = 0$). We have the estimates,
	\begin{align*}
		\norm{M}_{C^{k + 1, \al}(\Omega)}
			\le C \norm{\bomega}_{C^{k, \al}(\Omega)}, \quad
		\norm{M}_{H^{k + 1, p}(\Omega)}
			\le C \norm{\bomega}_{H^{k, p}(\Omega)},
	\end{align*}
	and
	$
		\sum_{j \ge 1} \abs{\la_j}
			\le C \norm{u}_{L^2(\Omega)}.
	$
	Moreover, if $\bomega$ is in the range of the curl then $\la_j = 0$
	for all $j$.
\end{theorem}
\begin{proof}
	All these observations follow from \cite{KMPT2000}, the explicit bound on $M$
	holding by the continuity of the solution
	operator $G^R$ defined in the proof of Corollary 3.2 of \cite{KMPT2000} and
	the comments in Section 5 of \cite{KMPT2000}.
	If $\bomega$ is in the range of the curl then the external fluxes vanish,
	which gives each $\la_j = 0$ as we can see in (2.2) of \cite{KMPT2000}
	(and see the comment immediately following the proof of Proposition 3.1
	of \cite{KMPT2000}).
\end{proof}

\begin{cor}\label{C:BSLawLikeBound}
	Assume that $\Gamma$ is $C^{k + 1, \al}$-regular, $k \ge 0$, and
	let $\bomega \in C^{k, \al}(\Omega)$
	(or $\bomega \in H^{k, p}(\Omega)$, $p \in(1, \iny)$)
	be in the range of the curl.
	Then there exists a unique
	$\uu \in H_0 \cap C^{k + 1, \al}(\Omega)$
	(or $\uu \in H_0 \cap H^{k + 1, p}(\Omega)$)
	for which $\curl \uu = \bomega$,
	and we have
	\begin{align*}
		\norm{\uu}_{C^{k + 1, \al}(\Omega)}
			\le C \norm{\bomega}_{C^{k, \al}(\Omega)}, \quad
		\norm{\uu}_{H^{k + 1, p}(\Omega)}
			\le C \norm{\bomega}_{H^{k, p}(\Omega)}.
	\end{align*}
\end{cor}
\begin{proof}
	Let $M$ be as in \cref{T:KMPT} and observe that $\dv M = \curl \vv$ for
	$\vv = (M_2^3, M_3^1, M_1^2)$.
	Solve
	\begin{align*}
		\begin{cases}
			\Delta p = \dv \vv
				&\text{in } \Omega, \\
			\grad p \cdot \n = \vv \cdot \n
				&\text{on } \Gamma,
		\end{cases}
	\end{align*}
	and let $\widetilde{\uu} = \vv - \grad p$.
	Then $\curl \widetilde{\uu} = \bomega$, $\dv \widetilde{\uu} = 0$,
	and $\widetilde{\uu} \cdot \n = 0$ on $\Gamma$. Moreover, $\dv \vv \in C^{k, \al}(\Omega)$
	and $\vv \cdot \n \in C^{k + 1, \al}(\Gamma)$, so elliptic estimates
	(as in item 3 of Lemma 2 in \cite{Koch2002})
	give $p \in C^{k + 2, \al}(\Omega)$. Letting $\uu = P_{H_0} \widetilde{\uu}$,
	and noting that $P_{H_0}$ is continuous in $C^{k + 1, \al}(\Omega)$
	by \cref{L:TamingHc}, we have
	\begin{align*}
		\norm{\uu}_{C^{k + 1, \al}(\Omega)}
			&\le C \norm{\vv}_{C^{k + 1, \al}(\Omega)} +
				C \norm{\grad p}_{C^{k + 1, \al}(\Omega)}
			\le C \norm{\vv}_{C^{k + 1, \al}(\Omega)}
			\le C \norm{\bomega}_{C^{k, \al}(\Omega)}
	\end{align*}
	by \cref{T:KMPT}. Similar estimates hold for Sobolev spaces.
\end{proof}

% \newpage
%\bibliography{Refs}
%\bibliographystyle{plain}

\def\cprime{$'$} \def\polhk#1{\setbox0=\hbox{#1}{\ooalign{\hidewidth
  \lower1.5ex\hbox{`}\hidewidth\crcr\unhbox0}}}

\end{document}